\documentclass[authoryear,preprint,11pt]{elsarticle}

\usepackage{setspace}
\usepackage{amssymb}
\usepackage{amsmath}
\usepackage{color}
\usepackage[dvipsnames]{xcolor}

\usepackage{eurosym}
\usepackage{lineno}
\usepackage{tabularx}
\usepackage{booktabs}
\usepackage{multirow}
\usepackage{float}
\usepackage{pdflscape}
\usepackage{hyperref}
\usepackage{afterpage}
\usepackage[left=3cm,right=3cm, top=4cm, bottom=4cm]{geometry}
\usepackage{geometry}

\usepackage[edges]{forest}
\usepackage{adjustbox}
\usepackage{tikz}
\usetikzlibrary{trees}
\newcommand{\Tau}{\mathrm{T}}
\newcommand{\cP}{\mathcal{P}}

\usepackage{titlesec}
\titleformat{\subsection}
  {\bfseries} 
  {\thesubsection} 
  {1em} 
  {} 

\titleformat{\paragraph}[runin]
  {\bfseries} 
  {\theparagraph} 
  {1em} 
  {} 
  [.] 

 \usepackage{lipsum}
 \makeatletter
 \def\ps@pprintTitle{%
 \let\@oddhead\@empty
 \let\@evenhead\@empty
\def\@oddfoot{}%
\let\@evenfoot\@oddfoot}

\usepackage{xcolor}
\usepackage{xargs}
\usepackage[colorinlistoftodos,prependcaption,textsize=tiny]{todonotes}

\usepackage{etoolbox}

\newcommand{\pubnotice}{%
  \par\noindent
   \setlength{\fboxsep}{8pt}%
  {\color{red}
  \fbox{%
    \parbox{0.95\textwidth}{%
      \normalfont\upshape
      This paper has been published as:\\
      Caselli, G., Iori, M., Ljubić, I.\ (2026).
      Bilevel optimization with sustainability perspective: A survey on applications.
      \textit{European Journal of Operational Research}, 332(2), 357--375.\\
      \href{https://doi.org/10.1016/j.ejor.2025.08.051}{https://doi.org/10.1016/j.ejor.2025.08.051}
    }%
  }%
}%
  \par\medskip
}

\makeatletter
\patchcmd{\pprintMaketitle}
  {\hrule\vskip12pt}
  {\pubnotice\hrule\vskip12pt}
  {}{}
\makeatother

\begin{document}

\begin{frontmatter}

\title{Bilevel optimization with sustainability perspective: A survey on applications}

\author[inst1,inst2]{Giulia Caselli\corref{cor1}}
\ead{giulia.caselli@unimore.it}
\author[inst1]{Manuel Iori}
\ead{manuel.iori@unimore.it}
\author[inst3]{Ivana Ljubić}
\ead{ivana.ljubic@essec.edu}
\cortext[cor1]{Corresponding author}

\affiliation[inst1]{organization={DISMI, University of Modena and Reggio Emilia},
            addressline={Via Amendola 2}, 
            city={Reggio Emilia},
            postcode={42122}, 
            country={Italy}}

\affiliation[inst2]{organization={DEMB, University of Modena and Reggio Emilia},
            addressline={Via Berengario 51}, 
            city={Modena},
            postcode={41121}, 
            country={Italy}}

\affiliation[inst3]{organization={IDS Department, ESSEC Business School},
            city={Cergy-Pontoise},
            country={France}}

\begin{abstract}
Bilevel optimization, a well-established field for modeling hierarchical decision-making problems, has recently intersected with sustainability studies and practices, resulting in a series of works focusing on bilevel optimization problems involving multiple decision makers with diverse economic, environmental, and social objectives.
This survey offers a comprehensive overview of sustainable bilevel optimization applications. 
First, we introduce the main concepts related to the nature of bilevel optimization problems 
and present some typical mathematical formulations for bilevel pricing problems that cover many of the collected applications.
Then, we review the most relevant works published in sustainable bilevel optimization, giving a classification based on the application domains and their association with well-known operations research problems, while briefly discussing the proposed solution methodologies. 
We survey applications on transportation and logistics, production planning and manufacturing, water, waste, and agriculture management, supply chains, and disaster prevention and response.
Finally, we outline a list of open questions and opportunities for future research in this domain.
\end{abstract}




\end{frontmatter}


\section{Introduction}
\label{sec:intro}
Hierarchical decision-making processes, involving interconnected and selfish players, often arise in the policy-making context
(e.g., organizational entities controlling local or private entities) or due to the inherent nature of decision processes (e.g., operational decisions sequencing in production). 
The prominent game theory concept of \textit{Stackelberg game} (see \citealt{stackelberg1934, stackelberg1952}), frequently used in economics to model hierarchical decision-making processes, has been introduced in mathematical optimization in the 70s (see \citealt{brackenmcgrill1973}), when the term \textit{bilevel optimization} is born.

Bilevel optimization allows to formalize decision processes where multiple decision makers are organized within a two-level hierarchy. 
Leader/follower decision problems permeate policy making in societal and industrial contexts, where sustainability concerns can shape the decisions of today for the survival of future generations. However, balancing economic prosperity with environmental concerns and social development is not a trivial task.
Operations research studies have recently progressed on problems and solutions for \emph{sustainable process optimization} (\citealt{JAEHN2016243}; \citealt{BARBOSAPOVOA2018399}), and bilevel optimization has been applied to sustainable practices since its first steps (see, e.g., \citealt{candlernorton1977}). Although the two topics are indeed connected, a clear picture on their formal interrelations is missing in the literature. 
Therefore, in this survey, we focus on major trends of bilevel optimization applications to problems related to sustainability, covering combinations of \textit{economic}, \textit{environmental}, and \textit{social} goals. 

\textit{Sustainability} is a broad and multidimensional concept that has evolved significantly over time. A big influence on the use of this term came from the 1983 United Nations Commission's on Environment and Development Brundtland Report (\citealt{Brundtland1987}) that defined sustainable development as ``development that meets the needs of the present without compromising the ability of future generations to meet their own needs''.
The United Nations' Sustainable Development Goals (SDGs) provide now a structured framework for sustainability across economic, environmental, and social dimensions (\citealt{UN2015}). 
The SDGs are 17 interrelated goals that address a range of global challenges, including poverty eradication, climate action, responsible consumption and production, and sustainable economic growth. 
These goals build upon the triple bottom line concept (\citealt{ elkington2002}), which defines sustainability through the three interconnected economic, environmental, and social dimensions.

Based on this principle, economic development should be pursued by companies and societies in a manner that benefits both the well-being of people and the health of the planet.
Sustainability-driven decision-making, particularly in policy and industrial settings, requires balancing economic interests with ecological and societal well-being. Bilevel optimization, by explicitly modeling hierarchical decision-making structures, provides a powerful tool for integrating these considerations. By incorporating sustainability principles into optimization models, researchers and decision makers can explore trade-offs between multiple objectives, addressing several SDGs, and facilitating solutions that support sustainable development.

Bilevel optimization problems follow a hierarchical structure where a leader, higher in the hierarchy, makes the first decision, and a follower responds by solving a nested optimization problem. The follower's problem is parameterized by a fixed leader decision, and, in turn, the lower-level solution affects either the feasible region or the objective function of the upper-level problem.
The follower’s problem must be solved to optimality to have a feasible solution for the bilevel problem. The global optimality of the lower-level problem reflects the \textit{rationality} of the follower's behavior, that is, the follower acts in favor of their objectives.
A well-posed bilevel formulation requires a non-empty lower-level feasible set, and the follower must select one among multiple optimal solutions once the leader decides (i.e., optimistic or pessimistic approach, discussed in Section \ref{sec:overview}). 
The solution methodologies for bilevel optimization problems strongly depend on the structure and properties of the lower-level problem. 
Much of the research since the 70s has focused on developing effective solution methods.
We refer the reader to \cite{KleinertLLS21} and \cite{beck2023survey} for recent surveys on exact approaches.
Real-world applications of bilevel optimization often involve complex lower-level problems and large-scale instances. Therefore, there has been an increasing interest in metaheuristic approaches over the past two decades. \cite{CAMACHOVALLEJO2024} present the first review of metaheuristic algorithms for bilevel problems, highlighting dominant methodologies and common practices. 
In this work, we provide a structured overview of current methodologies spanning reformulations, exact methods, and metaheuristics. 
Then, we focus specifically on real-world applications to derive a picture of notable bilevel problems from the sustainability perspective across different sectors to provide academics and practitioners with a screenshot of the state of the art on the use of bilevel programming in practice.

The application of bilevel optimization has extended across various fields such as the military industry, transportation, and energy management, as documented by \cite{Dempe2020} and \cite{sinhamalo2018}. 
\cite{Dempe2020} provides an exhaustive overview of topics of diverse bilevel applications, covering economics, energy and electricity, interdiction, and machine learning areas. 
Given space limitations and our focus on sustainability, we (i) exclude energy and electricity, intending to address them in a separate work due to the vast scope of the field, and (ii) select fields and top-journal studies specifically developed in sustainability terms, a focus not previously explored in bilevel optimization surveys. For a complete list of bilevel applications, we refer readers to \cite{Dempe2020}.
This study aims to offer a critical review of applications with sustainability perspective in bilevel optimization.

\begin{figure}
\footnotesize
\centering
\begin{adjustbox}{width=\linewidth}
\begin{tikzpicture}[
    applications/.style={
        text centered, text width=10em,
        text=black
    },
    sections/.style={
        text centered, text width=10em,
        text=black
    },
    subsections/.style={
        grow=down,
        xshift=-3.2em,  
        text width=10em,
        edge from parent path={(\tikzparentnode.205) |- (\tikzchildnode.mid west)}
    },
    level1/.style ={level distance=5em,anchor=west},
    level2/.style ={level distance=10em,anchor=west},
    level3/.style ={level distance=15em,anchor=west},
    level4/.style ={level distance=20em,anchor=west},
    level 1/.style={edge from parent fork down,sibling distance=15em,level distance=6em},
    level 2/.style={edge from parent fork down,sibling distance=10em}
]
    \node[anchor=south,applications](super){\textbf{Sustainable bilevel applications}}[]
    child{node [sections] {\textbf{Transportation and logistics}}
        child[subsections,level1] {node[anchor=west,level1] {Green location and transportation}}
        child[subsections,level2] {node[anchor=west,level2] {Sustainable routing}} 
    }
    child{node [sections] {\textbf{Production planning and manufacturing}}
        child[subsections,level1] {node[anchor=west,level2] {Carbon policies}}
        child[subsections,level2] {node[anchor=west,level2] {Green product design and manufacturing}}
    }
    child{node [sections] {\textbf{Water, waste, and agriculture management}}
        child[subsections,level1] {node[anchor=west,level2] {Water management}}
        child[subsections,level2] {node[anchor=west,level2] {Waste management}}
        child[subsections,level3] {node[anchor=west,level3] {Agriculture}} 
    }
    child{node [sections] {\textbf{Supply chains}}
        child[subsections,level1] {node[anchor=west,level1] {Green supply chain design and planning}}
        child[subsections,level2] {node[anchor=west,level2] {Eco-industrial parks}} 
    }
    child{node [sections] {\textbf{Disaster prevention and response}}
        child[subsections,level1] {node[anchor=west,level1] {Hazmat transportation}}
        child[subsections,level2] {node[anchor=west,level2] {Hazmat toll setting}} 
        child[subsections,level3] {node[anchor=west,level3] {Hazmat control policies}}
        child[subsections,level4] {node[anchor=west,level4] {Emergency planning and disaster response}}
    }
    ;
\end{tikzpicture}
\end{adjustbox}
\caption{Sustainable bilevel problems in various sectors.}
\label{fig:applications}
\end{figure}
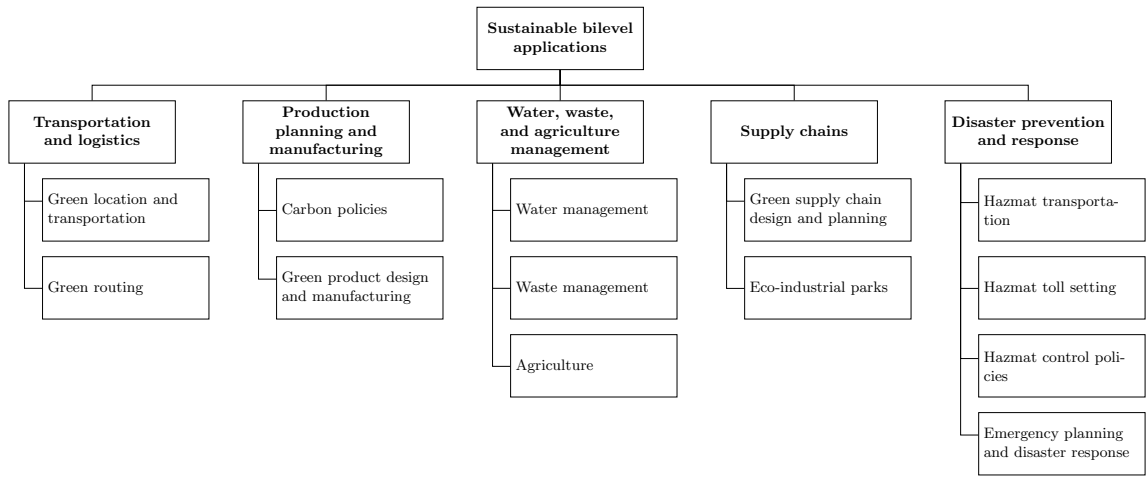

More in detail, this survey provides a taxonomy of applications with sustainability perspective  (see Figure \ref{fig:applications}) and analyzes the hierarchical nature of the bilevel game, identifying the players (leaders and followers) at each level and the competition or collaboration between them, especially in the form of an \textit{equilibrium}. It examines the players (e.g., policy makers, companies), their decisions, their goals by detecting whether the leader and follower problems are single- or multi-objective, which of the three dimensions of sustainability (namely, economic, environmental, and social) are considered in each level, which performance indicators are used to measure the objectives, and which solution methods are proposed for solving the bilevel problems.

Specifically, the study classifies sustainable bilevel applications according to their sector. 
We include works where either environmental or social goals are directly included in the players' objective functions (e.g., carbon emissions reduction, maximum equity in resource distribution), or where the implementation of the bilevel application brings positive impacts on the environment and society (e.g., green technology selection, natural disaster response).
\tikzstyle{every node}=[draw=black, thin, minimum height=4em]
\vspace{-.5cm}

We do not cover single-level optimization or other game theory concepts dealing with sustainability, nor bilevel applications without sustainable applications.
We focus on bilevel optimization applications within the operations research literature that apply mathematical programming or known metaheuristic frameworks, and exclude 
multi-level problems, i.e., optimization problems where players are in hierarchical games of three or more levels. 

The reminder of the paper is structured as follows. 
Section \ref{sec:overview} provides a brief overview on bilevel optimization, players, decisions, goals, and solution methods.
In Sections \ref{sec:transp}-\ref{sec:disaster}, we present the following classes of bilevel problems with relevant applications with sustainability perspectives: transportation and logistics, including a review of common mathematical formulations for bilevel network pricing problems (Section \ref{sec:transp}); production planning and manufacturing (Section \ref{sec:production}); water, waste, and agriculture management (Section \ref{sec:waste}); supply chains (Section \ref{sec:sc}); disaster prevention and response (Section \ref{sec:disaster}). 
Finally, Section \ref{sec:concl} identifies the research gaps in the existing literature and provides some guidance for future research.
A full list of acronyms is provided in \ref{app:acro} to improve the readability of the paper.

\section{Overview of bilevel applications}
\label{sec:overview}
\paragraph{Bilevel optimization}
A general optimistic bilevel optimization problem modeling a sequential game with two players (called the leader and the follower) is formulated as 
\begin{align}
\label{eq1} \min_{x \in X, y} & \quad F(x,y) \\
\label{eq2} \text{s.t.} & \quad G(x,y) \ge 0 \\
\label{eq3} & \quad y \in S(x), 
\end{align}
where $S(x)$ is the set of optimal solutions of the $x$-parameterized problem
\begin{align}
\label{eq4} \min_{y \in Y} & \quad f(x,y) \\
\label{eq5} \text{s.t.} & \quad g(x,y) \ge 0.
\end{align}
Models \eqref{eq1}-\eqref{eq3} and \eqref{eq4}-\eqref{eq5} define, respectively, the so-called upper-level (UL) and lower-level (LL) problems. The UL variables $x \in \mathbb{R}^{n_x}$ represent the decisions of the leader, and the LL variables $y \in \mathbb{R}^{n_y}$ represent the decisions of the follower.
The UL and LL objective functions are given by $F, f : \mathbb{R}^{n_x} \times \mathbb{R}^{n_y} \rightarrow \mathbb{R}$.
The constraint functions $G : \mathbb{R}^{n_x} \times \mathbb{R}^{n_y} \rightarrow \mathbb{R}^m $ and $g : \mathbb{R}^{n_x} \times \mathbb{R}^{n_y} \rightarrow \mathbb{R}^l $ map the sets of constraints of, respectively, the UL and LL problems. 
The sets $X \subseteq \mathbb{R}^{n_x}$ and $Y \subseteq \mathbb{R}^{n_y}$ represent the feasible decisions for the leader and the follower. 
The leader, positioned higher in the hierarchy, decides first anticipating the follower's decision outcome. The follower, lower in the hierarchy, reacts to the leader's decision afterwards. 
In applications studied in this paper, it is commonly assumed that the leader is optimistic, meaning that he/she controls the variables of the follower in such a way that, in case of multiple LL optimal solutions, the follower makes the best solution for the leader.
This is the most studied version of bilevel problems in our context of applications.
On the contrary, in the pessimistic version the leader assumes that the follower chooses the worst solution for the leader. 
While there are ongoing studies on the properties and solution methods for pessimistic problems (see, e.g., \citealt{liu2018pessimistic,liu2020pessimistic}), their applications fall primarily in risk-averse domains such as interdiction, principal-agent problems, and venture investment, which remain out of the scope of this work.

\begin{figure}
\footnotesize
\centering
\begin{adjustbox}{width=\linewidth}
\begin{tikzpicture}[
    games/.style={%
        text centered, text width=8em,
        text=black
    },
    types/.style={%
        grow=down,
        text centered, 
        text width=15em,
        text=black
    },
    objectives/.style={%
        grow=down,
        text centered, 
        text width=16em,
    },
    level 1/.style={edge from parent fork down,sibling distance=35em,level distance=9em},
    level 2/.style={edge from parent fork down,sibling distance=17em,level distance=9em}
]
    \node[anchor=south,games, level 1](super){\textbf{Bilevel games}}[]
     child[level distance=16.21em, edge from parent path={(\tikzparentnode.south) |- ++(0,-2.475em) -| (\tikzchildnode.north)}]{node[types] {\textbf{Leader/Follower (SLSF)} \\
        Central authority/Local authority \\
        Policy maker/Company \\
        Policy maker/Service consumer \\
        Manufacturer/Supplier \\
        SC1/SC2
    }}
    child{node[types] {\textbf{Multi-leader-follower \\ games}}
        child{node [types] {\textbf{Leader/Followers (SLMF)} \\
        Policy maker/Local authorities \\
        Policy maker/Companies
        }}
        child{node [types] {\textbf{Leaders/Follower (MLSF)} \\
        Local authorities/Policy maker \\
        Companies/Policy maker
        }}
        child{node [types] {\textbf{Leaders/Followers (MLMF)} \\
        Manufacturers/Transporters \\
        }} 
    }
    ;
\end{tikzpicture}
\end{adjustbox}
\caption{Taxonomy of bilevel games and players.}
\label{fig:games}
\end{figure}
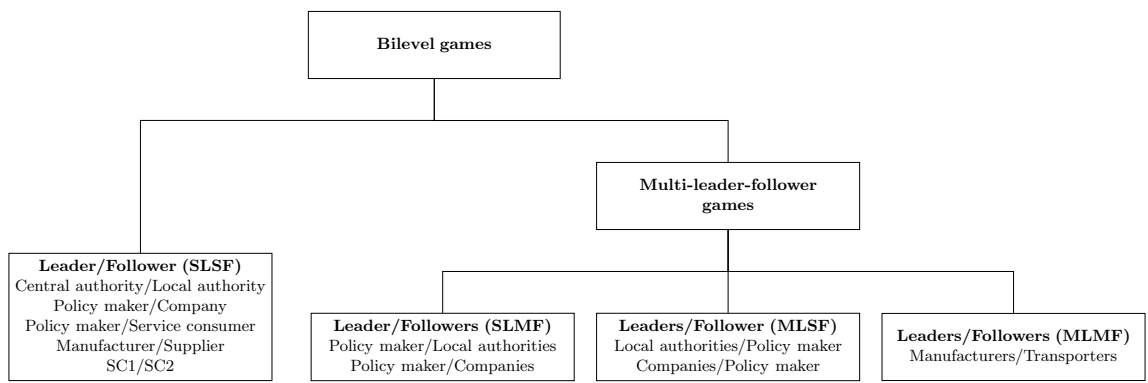

The class of single-leader-single-follower (SLSF) bilevel problems is well established (but still challenging) and represents several realistic situations involving one leader and one follower. 
In the so-called multi-leader-follower bilevel problems, multiple players are involved in the UL or LL problems. The classification includes single-leader-multi-follower (SLMF), multi-leader-single-follower (MLSF), and multi-leader-multi-follower (MLMF) problems (see Figure \ref{fig:games}).
In the intermediate setting of an SLMF game, there is a single leader but multiple followers, where either the followers are independent or they play a Nash game among themselves, given the leader's decision.
Let $y_j$ be the vector of decision variables of the $j$-th follower with $j=1, \dots, n_f$, where $n_f$ is the number of followers, 
and $x$ the leader's vector of variables.
A general SLMF game in the optimistic version is then 
\begin{align}
\label{eq1mf} \min_{x \in X, y_1, \dots, y_{n_f}} & \quad F(x,y) \\
\label{eq2mf} \text{s.t.} & \quad G(x,y) \ge 0 \\
\label{eq3mf} & \quad y:=(y_j)_{j= 1, \dots, n_f} \; \text{solves} \; \mbox{GNEP}(x), 
\end{align}
where $y$ is the vector of all decisions of all followers and GNEP$(x)$ is the set of \emph{generalized Nash equilibra} of the non-cooperative game among the $n_f$ followers, where the objective function and the feasible set of every follower explicitly depend on the decisions of the other followers. 
Indicating by $y_{-j}$ the decision variables of all followers except for those of the $j$-th follower, the LL problem of the $j$-th follower can be generally formulated as:
\begin{align}
\label{eq4mf} \min_{y_j} & \quad f_j(y_j,x,y_{-j}) \\
\label{eq5mf} \text{s.t.} & \quad y_j \in Y_j(x, y_{-j}), 
\end{align}
where  $f_j(y_j, x, y_{-j})$ and $Y_j(x, y_{-j})$ are the objective function and the feasible set of the $j$-th follower, the latter being defined by private constraints of the $j$-th follower and shared constraints of all followers. 
GNEPs arise quite naturally from standard Nash equilibrium problems if the players share some common
resource (e.g., a transportation link) or limitations (e.g., a common limit on the total pollution in a certain area).
The MLSF and MLMF games have more complex structures and are rarer in bilevel applications. We refer the reader to \cite{Aussel2020inbook} for general problem formulations and applications of MLSF and MLMF games.
In some applications, there are $n_f$ \emph{independent followers}, in which case the decision vector $y$ in model \eqref{eq1}-\eqref{eq3} is given as $y:=(y_j)_{j= 1, \dots, n_f}$. To  anticipate  the decisions of the followers, there are $n_f$ LL problems to be solved, i.e., the set $S(x)$ in \eqref{eq3} is the set of optimal solutions of $n_f$ $x$-parametrized problems:
\begin{align*}
 \min_{y_j \in Y_j}  \quad f_j(x,y_j),  
\quad \text{s.t.} \quad g_j(x,y_j) \ge 0,
\end{align*}
one for each $j \in \{1,\dots, n_f\}$, with objective function $f_j(x,y_j)$ and constraints set $g_j(x,y_j)$ mapping private constraints of the $j$-th follower only.

\paragraph{Players and decisions}
Figure \ref{fig:games} provides a schematic overview of bilevel games and the players involved.
In SLSF applications, the leader typically represents a central authority or a policy maker and the follower can be a local authority (e.g., the national authority sets an investment plan and the local authority designs the network layout), a company (e.g., the government sets taxes and subsidies and the company makes production and routing plans), or an average service consumer (e.g., the decision maker optimizes the network design and the user makes routing decisions). In other cases, the leader and the follower can be two similar entities hierarchically related and having conflicting goals, such as two members of the same supply chain (SC) or two SCs competing in the same market (``SC1'' and ``SC2''). In these cases, for instance, the former optimizes customer selection and routing decisions and the latter makes production planning decisions. 
Many real-world bilevel problems are modeled as SLMF games between one leader, a policy maker, and multiple followers, companies, local authorities, or end users. The policy maker sets prices, subsidies, budgets, taxes, allowances, or applies other policy tools to regulate some activities, typically concerning the regulation of market competition, carbon emissions, and social welfare. In consideration of such policies, companies make their production and distribution plans or local authorities design their tactics and operational policies or citizens, generally called end users, adopt their consumption behaviors and travel habits.
The MLSF and MLMF settings rarely appear in applications and in this survey we have selected very few of them.
Other bilevel applications of this kind are seen in the energy sector (see, e.g., \citealt{grimm2022tractable}).
\tikzstyle{every node}=[draw=black, thin, minimum height=4em, minimum width=15em]

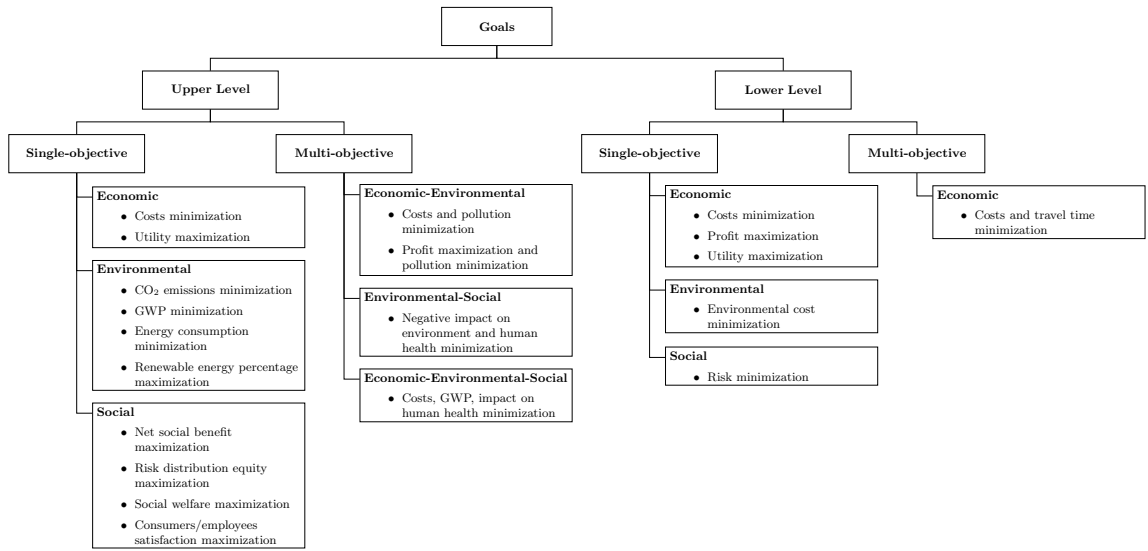
\begin{figure} 
\footnotesize
\centering
\begin{adjustbox}{width=\linewidth}
\begin{tikzpicture}[
    goals/.style={%
        text centered, text width=8em,
        text=black
    },
    levels/.style={%
        text centered, text width=10em,
        text=black
    },
    objectives/.style={%
        grow=down,
        xshift=.4cm,
        text width=16em,
        edge from parent path={(\tikzparentnode.south) |- (\tikzchildnode.mid west)}
    },
    level1/.style ={level distance=5em,anchor=west},
    level2/.style ={level distance=13.5em,anchor=west},
    level3/.style ={level distance=25.3em,anchor=west},
    level1mo/.style ={level distance=6em,anchor=west},
    level2mo/.style ={level distance=13.25em,anchor=west},
    level3mo/.style ={level distance=19em,anchor=west},
    level1r/.style ={level distance=5.7em,anchor=west},
    level2r/.style ={level distance=12em,anchor=west},
    level3r/.style ={level distance=16.7em,anchor=west},
    level1rmo/.style ={level distance=4.6em,anchor=west},
    level 1/.style={edge from parent fork down,sibling distance=45em,level distance=5em},
    level 2/.style={edge from parent fork down,sibling distance=21em}
]
    \node[anchor=south,goals](super){\textbf{Goals}}[]
    child{node [levels] {\textbf{Upper Level}}
        child{node [levels] {\textbf{Single-objective}}
            child[objectives,level1] {node[anchor=west,level2] {\textbf{Economic}
            \begin{itemize}
                \item Costs minimization
                \item Utility maximization
            \end{itemize}
            }}
            child[objectives,level2] {node[anchor=west,level2] {\textbf{Environmental}
            \begin{itemize}
                \item CO\textsubscript{2} emissions minimization
                \item GWP minimization
                \item Energy consumption minimization
                \item Renewable energy percentage maximization
            \end{itemize}
            }}
            child[objectives,level3] {node {\textbf{Social}
            \begin{itemize}
                \item Net social benefit maximization
                \item Risk distribution equity maximization
                \item Social welfare maximization
                \item Consumers/employees satisfaction maximization
            \end{itemize}
            }}
            }
        child{node [levels] {\textbf{Multi-objective}}
            child[objectives,level1mo] {node {\textbf{Economic-Environmental}
            \begin{itemize}
                \item Costs and pollution minimization
                \item Profit maximization and pollution minimization
            \end{itemize}
            }}
            child[objectives,level2mo] {node {\textbf{Environmental-Social}
            \begin{itemize}
                \item Negative impact on environment and human health minimization
            \end{itemize}
            }}
            child[objectives,level3mo] {node {\textbf{Economic-Environmental-Social}
            \begin{itemize}
                \item Costs, GWP, impact on human health minimization
            \end{itemize}
            }}}  
    }
    child{node [levels] {\textbf{Lower Level}}
        child{node [levels] {\textbf{Single-objective}}
            child[objectives,level1r] {node {\textbf{Economic}
            \begin{itemize}
                \item Costs minimization
                \item Profit maximization
                \item Utility maximization
            \end{itemize}
            }}
            child[objectives,level2r] {node {\textbf{Environmental}
            \begin{itemize}
                \item Environmental cost minimization
            \end{itemize}
            }}
            child[objectives,level3r] {node {\textbf{Social}
            \begin{itemize}
                \item Risk minimization
            \end{itemize}
            }}}
        child{node [levels] {\textbf{Multi-objective}}
            child[objectives,level1rmo] {node {\textbf{Economic}
                        \begin{itemize}
                \item Costs and travel time minimization
            \end{itemize}
            }}}  
    };
\end{tikzpicture}
\end{adjustbox}
\caption{Taxonomy of economic, environmental, and social objective functions in bilevel applications.}
\label{fig:goals}
\end{figure}

\paragraph{Goals}
We analyze the three dimensions of sustainability (namely, economic, environmental, and social) within the UL and LL problems of the selected bilevel applications to identify common roles and behaviours of the game players. Figure \ref{fig:goals} provides a schematic overview.
The LL problems typically have one economic goal, such as costs minimization, profit maximization, or some utility function maximization for the followers. Sometimes the LL problem is multi-objective with two utility goals such as cost and travel time minimization in network optimization problems. 
Rarely, when the follower is a local authority relating to a central authority in the UL for the design of SCs or industrial parks, an environmental goal is included in the LL problem, such as pollution cost within the total cost minimization or the solely environmental cost minimization.
In specific applications with a similar setting concerning, for instance, the hazardous material transportation,
a social goal is included in the LL when the follower is again a local area authority relating to a central authority in the UL for the design of the transportation network. The follower aims to minimize the total risk of hazmat transportation on the population of the area (e.g., increased accident probabilities), while the central authority maximizes the risk equity across areas.
On the contrary, the UL problems show a variety of single- and multi-objective functions including the three dimensions of sustainability. Very often, the leader is a policy maker whose aim is to improve the economic, environmental, and social conditions of the system (e.g., a country or a region, an industrial sector, a public service). The single-level objective functions are economic (total costs minimization, total utility maximization), environmental (minimization of CO\textsubscript{2} emissions, global warming potential (GWP), energy consumption, and maximization of renewable energy percentage), or social (maximization of net social benefit, risk distribution equity, social welfare, consumer or employee satisfaction). Economic and environmental or environmental and social goals are often combined in bi-objective UL problems. For instance, impact on land and human health are minimized in agriculture bilevel applications. In general, pollution minimization is combined with costs minimization or profit maximization in public network design, SC management, and carbon policy setting. In some cases, the three dimensions are all included in multi-objective UL problems.
\tikzstyle{every node}=[draw=black, thin, minimum height=3em]

\normalsize
\paragraph{Solution approaches} 
Bilevel optimization problems are notoriously hard to solve. For instance, \cite{hansen1992new} showed that bilevel problems where UL and LL are linear programming (LP) problems (i.e., bilevel LPs) are strongly $\mathcal{NP}$-hard. 
For bilevel MILPs, even checking the feasibility of a given point is  $\mathcal{NP}$-hard (\citealt{KleinertLLS21}).
Solution methods for bilevel problems strongly depend on the structure and properties of the LL problem (e.g., continuous linear or convex, mixed-integer linear, non convex) and on the coupling between the UL and LL problems (e.g., optimistic or pessimistic version).

\textit{Single-level reformulations.}
Assuming convexity and a suitable constraint qualification of the LL problem for all possible UL decisions, an SLSF problem is usually reformulated into a single-level problem, that is, a single-level reformulation (SLR) by imposing Karush-Kuhn-Tucker (KKT) conditions, strong duality conditions for the follower's solution $y$,
or variational inequalities.  
In the first case, we will refer to the KKT reformulation (KKTR). 
Solving the SLR to global optimality guarantees finding an optimistic globally optimal bilevel solution.
Another standard reformulation is the lower-level value function reformulation (LLVFR), which relies on the LL optimal value function definition (\citealt{outrata1990}).
For instance, \cite{Fischer2022} recently proposed a Newton method for finding stationary points of optimistic bilevel problems based on the LLVFR with successful results on well-known instances.
Notably, the LLVFR does not require the convexity and constraint qualification assumptions for equivalence between the SLR and the bilevel problem optimality. 
In general, the value function is not differentiable resulting in a non-smooth optimization problem.
The LLVFR guarantees both local and global equivalence to the bilevel problem (\citealt{zemkoho2021}), whereas KKTR only ensures global equivalence. Indeed, locally optimal solutions of the KKTR are not necessarily locally optimal for the original bilevel problem due to the introduction of implicit variables in the reformulation (\citealt{dempe2012bilevel}).
To the best of our knowledge, no definitive preference has been established between KKTR and LLVFR, and both approaches have been applied in the literature. Other duality-based reformulations have been analyzed in \cite{dempe2024duality}.
There are various techniques for solving such obtained non-convex nonlinear models (i.e., SLRs) to optimality.

\textit{Exact methods.}
Exact solution approaches applied to LP, MILP, and mixed-integer nonlinear programming (MINLP) bilevel problems include Branch-and-Bound (B\&B), Branch-and-Cut (B\&C), Column Generation (CG), Cutting Plane (CP), or Benders' Decomposition (BD) algorithms.
Typically, a suitable SLR is derived, to which bilevel-specific B\&B or B\&C techniques are applied. Some duality-based SLRs can also be solved directly using general-purpose solvers.
The penalty function method and approximation heuristics are typically applied to the KKTR. Enumeration algorithms, including the well-known \textit{$k$-th best} algorithm, are sometimes applied to bilevel linear programs.
For bilevel MILPs, the common approach is to first solve the high-point relaxation and then iteratively discard bilevel infeasible solutions using adapted branching and cutting techniques, see \cite{FischettiLMS18,TahRalDeN20}. Benders’ decomposition can be adapted for bilevel MILPs, assuming high-point relaxation compactness (\citealt{BolRal22}).
For non-convex bilevel MILPs, a few decomposition and branching approaches have been proposed, see, e.g., \cite{GaarLLST24} and further references therein. 
\\
We refer to the surveys by \cite{KleinertLLS21} for a complete overview of exact approaches for bilevel optimization.
For bilevel optimization under uncertainty, we refer the reader to \cite{beck2023survey}. In this context, both the stochastic optimization (SO) and robust optimization (RO) approaches can be used to model uncertainty at each of the two decision levels. 

\textit{Metaheuristics.}
In most applications dealing with real-world bilevel problems, LL problems are non-convex or give intractable SLRs. Therefore, metaheuristics like genetic and evolutionary algorithms (GA; EA), (multi-objective) particle swarm optimization (PSO; MOPSO), and the non-dominated sorting genetic algorithm (NSGA-II) are often used and are especially effective in the case of multi-objective UL or LL problems. Hybrid algorithms are also implemented by combining heuristics and LP/MILP models.
Common metaheuristic frameworks are: (i) for each UL decision, the LL problem is solved either exactly or heuristically (i.e.,  ``nested approach''); (ii) an SLR is derived and then solved by a metaheuristic algorithm (most commonly genetic ones). 
Less common are the co-evolutionary (or interactive) framework, where information is exchanged iteratively between the UL and LL levels, and the transformation of the bilevel problem into a bi-objective one. 
In general, metaheuristic methods focus on obtaining high-quality LL solutions rather than solving the LL problem to optimality.
For multi-objective bilevel problems, the vector optimization problem is often replaced with a scalar-valued optimization problem (scalarization or ``weighted sum'' method). Some results on (weak) efficiency of solutions are provided in \cite{CAMACHOVALLEJO2024}, while genetic and evolutionary methods remain dominant in applications.   
Their survey also provides an in-depth overview of metaheuristic methods for bilevel optimization.
While these methods enhance scalability and computational feasibility, they may lead to solutions that violate follower’s rationality condition. Great research efforts are addressed to explore equivalent single-level reformulation and exact methods. However, such approaches face substantial challenges for problems with both simple and complex structures. Balancing the strengths of metaheuristics with the rigor of exact methods remains an important open challenge in the field.

\section{Transportation and logistics}
\label{sec:transp}
Network design is a major class of combinatorial optimization problems conceptually represented as the selection of a subset of links in a graph, with numerous application developments in transportation and logistics (see \citealt{Cordeau2021}).
A substantial body of bilevel optimization research has focused on problems with a network structure (see, e.g., \citealt{brotcorne2008montreal, Labbe2021}), with most applications, as highlighted in this work, following this framework.
In this section, we present three typical mathematical formulations for bilevel network pricing problems focusing on the optimistic setting. 
These formulations cover the majority of applications 
discussed here, including transportation, facility location, and logistics. They are also widely used in hazardous material toll-setting problems (see Section \ref{subsec:haz-toll}).
We then explore relevant applications in the field, highlighting how each study adapts the general models to specific contexts.
In this work, we do not survey non-sustainable applications of bilevel network design.
The works we selected address bilevel network pricing and design applications in sustainable transportation and logistics. 
These works align with SDG 9 (Industry, Innovation, and Infrastructure) by promoting resilient infrastructures, sustainable manufacturing, and CO\textsubscript{2} emissions reduction, and SDG 11 (Sustainable Cities and Communities) by addressing urban disaster risk, enhancing public transport, reducing air pollution, and optimizing road networks.
We include studies related to ``green'' extensions of problems commonly encountered in this domain, like the Facility Location Problem (FLP) and the Vehicle Routing Problem (VRP).
Governmental policies for reducing the environmental impact of, for instance, industrial production and supply chains, and improving social welfare are modeled as pricing problems.
Network pricing decisions are crucial for airline, freight, and urban transportation, as well as telecommunication and service industries, especially in modern markets where intense price competition from economic deregulation on one hand, and environmental and social concerns from sustainable development regulations and practices on the other hand, coexist with network modifications.

\subsection{Typical network pricing formulations}
\label{subsec:models}
Our goal is to analyze the emerging bilevel problems with real-world applications that relate to sustainability practices and present general mathematical formulations for the basic problems. 
The given formulations can 
serve as a reference for researchers dealing with specific problems in the class and can be used as a starting point for further generalizations.

\subsubsection{Price setting problems}
\label{subsec:tollmodel}
The first model is a basic linear price setting problem with one leader and one or multiple independent followers. 
Let $T \in \mathbb{R}^{n_1}$ be the price vector, representing tax, tariff, or subsidy (if negative) decisions of the leader and $\Tau$ be a polyhedron representing the set of feasible tax decisions for the leader. Let  $y \in \mathbb{R}^{n_1}$ and $z \in \mathbb{R}^{n_2}$ be the decision variables of the follower on taxed and untaxed activities, respectively. In addition, let $c_1 \in \mathbb{R}^{n_1}$ and $c_2 \in \mathbb{R}^{n_2}$ be the LL objective function coefficients associated to activities $y$ and $z$, respectively.
The price setting SLSF problem is formulated as 
\begin{equation}
\label{eq8} \hspace{-3cm} \max_{T \in \Tau, y}  \quad Ty 
\end{equation}
\vspace{-.3cm}
\begin{subequations}
\begin{equation}
\label{eq9} \hspace{2cm} \text{s.t.}  \quad (y,z) \in \text{arg} \min_{y,z}  \quad (c_1 + T)y + c_2z \\
\end{equation}
\begin{equation}
\label{eq10}  
 \hspace{1cm} \text{s.t.}  \quad g(y,z) \le 0  \\
\end{equation}
\begin{equation}
\label{eq11} \hspace{.4cm}
  \hspace{1.5cm} \quad y \in Y, z \in Z
\end{equation}
\end{subequations}
where the bilinear objective function of the leader, given by \eqref{eq8},  maximizes tax revenues obtained from taxable activities of the follower, priced according to the tax values of $T$. The LL problem \eqref{eq9}-\eqref{eq11} typically models a cost-minimization problem (such as shortest path, or routing of multiple commodities). Constraints \eqref{eq10}-\eqref{eq11} model the feasible space of the follower (observe that the feasible LL set is not affected by the decisions of the leader). The interaction between the leader and the follower is restricted to the cost coefficients in the objective function of the LL \eqref{eq9}. If the tax $T$ is set too high, the follower will prefer to avoid taxed activities, which in return will diminish the overall revenue of the leader (see \citealt{LabbeViolin:2016} for more details). Constraints \eqref{eq11} are typically given as polyhedra with possible integer restrictions on some of the decision variables.
For bilevel problems of this type (the follower's constraint set does not depend on the leader's decision and all functions are linear) \cite{DEMPE2014solution} propose an effective solution method based on the LLVFR (optimal value function reformulation).
The model can be easily generalized to multiple independent followers (cf.\ Section \ref{sec:overview}). Moreover, the linear dependence between the LL activities and the prices set by the leader can also be generalized so that the objective function of the follower is given as 
\[ \min C_1(y,T) + C_2(z) \]
where $C_1$ and $C_2$ are nonlinear cost functions, the former representing the cost for the follower's taxed activities and the latter the cost of the non-taxed activities.

\subsubsection{Joint pricing and network design}
In the joint network design and pricing problem (see, e.g.,  \citealt{Brotcorne2008}), 
the leader seeks to design an infrastructure (typically modeled as a network) and to decide on the prices for the follower's activities if this infrastructure is used. 
The leader decides on the tariff as in pricing model \eqref{eq8}-\eqref{eq11} (vector $T$ of variables), and on the design of the network. The latter decision is represented by the binary variables vector $x \in X$, where set $X$ contains budget constraints and other restrictions on the feasible infrastructure. If the follower uses the infrastructure (modeled using $y$ variables at the LL), the activity will be taxed accordingly. The follower can also opt-out for non-taxed activities (modeled by $z$ variables as above). 
The joint pricing and network design SLSF problem 
is then
\begin{equation}
\label{eq16} \hspace{-3cm} \max_{T \in \Tau, \; x \in X, y}  
\quad Ty 
\end{equation}
\vspace{-.5cm}
\begin{subequations}
\begin{align}
\label{eq17} \hspace{2cm} \text{s.t.}  \quad (y,z) \in \text{arg} \min_{y,z}  \quad C_1(y,T) + C_2(z) \\
\label{eq18}  
 \hspace{1cm} \text{s.t.}  \quad g(y, z) \le 0 \\
\label{eq19}  
 \hspace{2cm} \quad g_1(x, y) \le 0 \\
\label{eq20} \hspace{1.5cm} \quad y \in Y, z \in Z.
\end{align}
\end{subequations}
The LL problem \eqref{eq17}-\eqref{eq20} is parameterized by the UL tariff variables $T$ and the UL variables $x$, which determine the subnetwork on which the LL problem is solved.
As above, \eqref{eq17} models the cost function of the follower and constraints \eqref{eq18} guarantee the feasibility of the follower's problem by performing taxed or non-taxed activities. The major difference in the LL problem with respect to the pure price setting model is in constraints \eqref{eq19}: they ensure that the follower can use the infrastructure only if it is installed by the leader. Typically, these constraints on network arcs are given as $y_{a} \le U_a\cdot x_a $, where $a$ represents an arc in the network, and $U_a$ is the maximum arc capacity. These constraints are linking network design decisions of the leader and the (routing) decisions of the follower (i.e., the follower can use the network only after the leader has undertaken design decisions). Also this problem can be easily generalized to the setting with multiple independent followers. 

\subsubsection{Pricing with equilibrium constraints}
\label{subsec:eqmodel}
Bilevel pricing problems with a single leader and multiple followers sometimes contain the so-called network equilibrium problem in the LL.
We illustrate this type of problem as an example of a transportation problem in which the LL problem is modeled as a Wardrop equilibrium. 
Departing from the price setting model given by \eqref{eq8}-\eqref{eq11}, the leader still has 
to decide on the prices of the arcs in the network (decisions $T\in \Tau$, as above). However, followers are no longer independent players but interact with each other so that they are in a simultaneous equilibrium (i.e., no player can unilaterally improve their objective by changing their decision, see \citealt{wardrop1952road}).
Each follower represents a set of drivers sharing the same origin and destination.
Typical applications arise in traffic management of congested urban areas where a Wardrop equilibrium 
aims at finding the minimum travel time (or cost) of all followers.
This is because each driver decides to use the least duration path and therefore, at equilibrium, all routes with the same origin and destination have the same traveling time.  The deterministic user equilibrium assumes that the followers have complete and accurate information about the available paths and the network flows are stable over time. 
In the presence of incomplete information, we deal with stochastic models (see \citealt{patriksson2015traffic}).

Following the notation of a general equilibrium model from \cite{Morandi:2024} and bilevel network design models from \cite{Labbe2021}, we now define the LL Wardrop equilibrium problem which is parametrized in the pricing UL variables $T$. Let the graph $G = (V, A)$ represent the network, where $A$ is the set of arcs, and let $C$ be the set of commodities, each representing a given origin-destination (OD) pair with demand $d_c > 0$ (multiple drivers may share the same OD pair). 
For each commodity $c \in C$, let  $\cP_c$ be the set of all paths connecting their origin with their destination, and $p \in \cP_c$ be the index of a path connecting the origin and destination nodes of commodity $c$. 
In the LL model stated below, variables $y_a$ represent the aggregated flow on each arc  $a \in A$ of the network, while $z_p$ represent the flow of the OD pair $c$ on path $p \in \cP_c$.  
The average travel time on the network is modeled by a link performance function $f_a(y_a, T)$ (separable function), which depends on the flow on each arc and the leader's decision on arc prices. The travel time function also includes the so-called practical capacity parameters, introduced to measure flow congestion on each arc.
As in previous models, vector $T$ indicates UL variables for pricing decisions.

The LL results in a user equilibrium model (parametrized in $T$), also called \textit{traffic assignment problem} (see, e.g., \citealt{dempezemkoho2012}), formulated as
\begin{subequations}\label{equilb}
\begin{align}
\label{eq17eq} & \min_{y, z} \quad \sum_{a \in A} \int_0^{y_a}f_a(\omega,T)d\omega \\
\label{eq18eq} & \text{s.t.} \quad \quad y_a = \sum_{c \in C}\sum_{\substack{p \in \cP_c: \\ a \in p}} z_p \quad && \forall a \in A\\
\label{eq19eq} & \quad \quad \quad  d_c = \sum_{p \in \cP_c} z_p \quad && \forall c \in C\\
\label{eq22eq} & \quad \quad \quad  z_p \ge 0 \quad && \forall c \in C, p \in \cP_c.
\end{align}
\end{subequations}

The model considers each OD pair as a commodity $c \in C$. 
At equilibrium, all used routes from an origin to a destination are characterized by the same average travel time. This equilibrium corresponds to a Nash equilibrium in a game with a large amount of players (i.e., followers). Traffic tends to settle down towards equilibrium where no travelers may want to change route. 
The KKT conditions of model \eqref{equilb} correspond to Wardrop equilibrium conditions. 
The objective function \eqref{eq17eq} of the equilibrium model
integrates the travel time function for each follower over each arc flow (i.e., the higher the congestion, the longer the travel time).  
It is the sum over all arcs of the integral between
0 and the arc flow $y_a$ of the latency function $f_a$. The function  $f_a$ is assumed to be continuous, positive, and increasing with flow $y_a$ to guarantee solution existence and the avoidance of cycles, and depends on the pricing decision $T$ of the leader.
Constraints \eqref{eq18eq} determine the arc flows $y_a$ as the sum of the flows of all paths passing through the arc $a \in A$. Constraints \eqref{eq19eq} guarantee demand $d_c$ satisfaction for each commodity $c \in C$.
Constraints \eqref{eq22eq} state the domain of the LL variables.
Common solution approaches to bilevel problems where the LL is modeled as a traffic assignment problem include sensitivity analysis, KKTR, and LLVFR (\citealt{dempezemkoho2012}).
A more general approach for modeling Wardrop equilibrium using mathematical programming with equilibrium constraints can be found in \cite{brotcorne2008montreal}.

\subsection{Green location and transportation}
\label{subsec:loc_and_transp}
In transportation bilevel problems, it is common that the government or network operator acts as the leader making decisions related to network design and price setting, while the users of the network act as the followers, making individual route choices.
Given the growing importance of sustainability in transportation, recent bilevel problems in this field deal with sustainability, either directly (by explicitly minimizing carbon emissions and promoting green mobility) or indirectly (through policies that encourage the adoption of sustainable transport modes).
Sustainable transportation, in our context, refers to alternative transport modes that aim to reduce carbon emissions, such as bicycles, electric vehicles, and efficient public transport systems.
The selected papers focus on bilevel problems that support sustainable transportation by addressing infrastructure investments, pricing mechanisms, and regulatory constraints that promote transportation at low emissions or with positive social impact.
Table \ref{tab:location} provides a classification of the referenced works based on their application, the type of bilevel game and whether the LL problem is modeled as an equilibrium problem (``Eq.''), players and objectives in the UL and LL problems, and the proposed solution method.
\afterpage{
\clearpage
\newgeometry{margin=3cm}
\begin{landscape}
\begin{table}[ht]
    \tiny
    \centering
    \caption{Summary of papers on bilevel applications focused on green bilevel location and transportation.
    }
    \label{tab:location}
    \begin{tabularx}{\linewidth}{p{0.10\textwidth}p{0.17\textwidth}p{0.05\textwidth}p{0.02\textwidth}p{0.15\textwidth}p{0.15\textwidth}p{0.22\textwidth}p{0.22\textwidth}X}
        \toprule
        Application  & References & \multicolumn{2}{c}{Bilevel game} & \multicolumn{2}{c}{Players} & \multicolumn{2}{c}{Objectives} & Method\\
        \cmidrule(lr){3-4}\cmidrule(lr){5-6}\cmidrule(lr){7-8} 
        & & Class & Eq. & UL & LL & UL & LL \\
        \cmidrule(lr){1-1} \cmidrule(lr){2-2}\cmidrule(lr){3-4}\cmidrule(lr){5-6}\cmidrule(lr){7-8}   \cmidrule(lr){9-9}        
        \multirow{5}{*}{\parbox{0.10\textwidth}{Highways network design and pricing}} &\cite{Ben-AyedBBL92} & SLSF & & Government & Average road user & Min total cost & Max utility & HEUR  \\
         & \cite{SinhaMD15} & SLSF/ SLMF & & Government & Average road user/Users & Max revenue; Min pollution & Min total cost; Min total travel time & MOEA \\
        &\cite{wang2014bilevel} & SLMF & \checkmark & Tolls DM & Road users & Min travel time; Min emissions; Min impact on health & Min travel cost; Min travel time & NSGA-II \\        
        & \cite{wen2016minimizing} & SLMF & \checkmark & Tolls DM & Road users & Min emissions; Min travel time & Min total cost & SA \\    
        & \cite{zhao2016bi} & SLMF & \checkmark & Tolls DM & Road users & Min emissions & Min total travel time & EA+AON \\
        \cmidrule(lr){1-1} \cmidrule(lr){2-2}\cmidrule(lr){3-4}\cmidrule(lr){5-6}\cmidrule(lr){7-8}\cmidrule(lr){9-9}
         \multirow{5}{*}{\parbox{0.10\textwidth}{Bike lanes network design}} & \cite{sohn2011multi} & SLMF & \checkmark & Network designer & Motorists and cyclists(or government) & Min motorists time and min cyclists time(or max motorist speed) & Min travel time & NSGA-II\\
         & \cite{bagloee2016bicycle} & SLMF & \checkmark & Network designer & Motorists and cyclists & Min total travel time &  Min travel time & B\&B\\ 
          & \cite{gaspar2015bilevel}  & SLMF & \checkmark & Network designer & Motorists, cyclists, and bus users & Max number of cyclists & Min travel time & HEUR \\
         & \cite{LiuSJ22} & SLMF & \checkmark & Network designer & Cyclists & Max total utility &  Max utility & KKTR \\
         & \cite{rashidi2016optimal} & SLMF & \checkmark & City planner & Motorists, public transit passengers, pedestrians & Min total travel cost (incl. safety) & Min travel cost (incl. safey) & HEUR \\
        \cmidrule(lr){1-1} \cmidrule(lr){2-2}\cmidrule(lr){3-4}\cmidrule(lr){5-6}\cmidrule(lr){7-8}\cmidrule(lr){9-9}
         \multirow{3}{*}{\parbox{0.10\textwidth}{Green facility location}} & \cite{ChenQMU20} & SLMF  & \checkmark & EV service provider & EV drivers & Min construction cost and min total drivers time & Min time & KKTR+IA \\ 
        & \cite{zhou2023location} & SLMF & \checkmark & Network designer & Users & Max transfer flow; Max transit priority &  Min travel time & HEUR \\
         & \cite{GangTLXSY15} & SLMF &  & Government & Enterprises & Min pollution; Min total cost & Min total cost & PSO \\
        \cmidrule(lr){1-1} \cmidrule(lr){2-2}\cmidrule(lr){3-4}\cmidrule(lr){5-6}\cmidrule(lr){7-8}\cmidrule(lr){9-9} 
         \parbox{0.10\textwidth}{Air passenger transportation} & \cite{qiu2020carbon} & SLMF &  & Goverment & Airlines & Max net social benefit &  Max profit & GA \\
        \bottomrule
    \end{tabularx}
\end{table}
\begin{table}[H]
    \tiny
    \centering
    \caption{Summary of papers on bilevel applications focused on sustainable routing. 
    }
    \label{tab:routing}
    \begin{tabularx}{\linewidth}{p{0.10\textwidth}p{0.17\textwidth}p{0.05\textwidth}p{0.02\textwidth}p{0.15\textwidth}p{0.15\textwidth}p{0.22\textwidth}p{0.22\textwidth}X}
    \toprule
    Application & References & \multicolumn{2}{c}{Bilevel game} & \multicolumn{2}{c}{Players} & \multicolumn{2}{c}{Objectives} & Method\\
    \cmidrule(lr){3-4}\cmidrule(lr){5-6}\cmidrule(lr){7-8}
    & & Class & Eq. & UL & LL & UL & LL &  \\
    \cmidrule(lr){1-1}\cmidrule(lr){2-2}\cmidrule(lr){3-4}\cmidrule(lr){5-6}\cmidrule(lr){7-8} \cmidrule(lr){9-9}       
    \multirow{2}{*}{\parbox{0.10\textwidth}{Pollution Routing Problem}} & \cite{nath2018novel} & SLMF & & Depot & Vehicles & Min fuel consumption; Min total travel distance & Min travel distance & NSGA-II+GA \\
    & \cite{qiu2020carbonpricing} & SLSF & & Authority & Company & Min CO\textsubscript{2} emissions & Min total (fuel, emissions, operating) costs & HEUR \\
    \cmidrule(lr){1-1}\cmidrule(lr){2-2}\cmidrule(lr){3-4}\cmidrule(lr){5-6}\cmidrule(lr){7-8} \cmidrule(lr){9-9} 
    \multirow{2}{*}{\parbox{0.10\textwidth}{School bus routing}} 
    & \cite{PARVASI2017}  & SLMF &  & Public transportation company & Students & Max profit & Min cost & Hybrid Alg \\ 
    & \cite{calvetegale2023}  & SLMF &  & School & Students & Min total cost & Max preference & Dual. SLR; HEUR \\ 
    \bottomrule
    \end{tabularx}
\end{table}
\end{landscape}
\restoregeometry
}

\paragraph{Highway network design and pricing}
\label{par:hightoll}

\cite{Ben-AyedBBL92} and \cite{SinhaMD15} introduce transportation bilevel problems focused on highway pricing and network design. These problems fall under the SLSF class, where the government acts as the leader and the average user serves as the follower. 
\cite{Ben-AyedBBL92} present early bilevel applications in pure network design, addressing rural highway systems in developing countries. The government aims to increase the capacity of the existing highway network, while the follower solves a flow optimization problem.
This early application aligns with sustainable development and the SDGs by promoting equitable access to transportation in rural areas, improving infrastructure resilience, and minimizing environmental impact.
The authors develop an iterative heuristic algorithm and solve a realistic instance from a case study in Tunisia.
\cite{SinhaMD15} extend the pricing model \eqref{eq8}-\eqref{eq11} from Section \ref{subsec:tollmodel} for the highway pricing problem with multi-objective UL and LL. In particular, the leader aims to maximize revenues ($F_1$, defined in \eqref{eq8}) while minimizing pollution costs ($F_2(y,z)=a_1y+b_1z$). At the LL, the follower aims to minimize total travel costs ($f_1$, defined in \eqref{eq9}) and total travel time ($f_2(y,z)=a_2y+b_2z$). 
A Pareto-optimal front of decisions balancing revenues and pollution is given for the UL problem, and an LL Pareto-optimal front balancing travel cost and time is given for any UL optimal decision.
The authors propose a multi-objective evolutionary algorithm (MOEA) and solve an illustrative example. They also extend the problem to the SLMF version.

In the presence of multiple followers (SLMF), sustainable tolls setting problems are often modeled as bilevel pricing problems with equilibrium (see model \eqref{eq17eq}-\eqref{eq22eq} of Section \ref{subsec:eqmodel}).
These models are frequently extended to incorporate multiple objectives in the UL, focusing not (only) on the leader's profit maximization, but on various performance measures of economic, environmental, and health sustainability in transportation.
In the UL, the leader is responsible for setting tolls.
The followers are the travelers of the system and the LL models an equilibrium with minimum travel time/cost for every traveler.
\cite{wang2014bilevel} study the bilevel tolls setting problem with three minimization objectives in the UL (namely, total system travel time, CO\textsubscript{2} emissions produced by fuel consumption, and human pollutant uptake), and a bi-objective equilibrium in the LL. The bi-objective equilibrium model is obtained by expressing the objective function $f_a(y_a, T)$ in equation \eqref{eq17eq} as the sum of two separable functions $T_a(y_a)$ and $M_a(T)$, respectively, for the total travel time and total toll cost, and exploiting results from bi-objective Wardrop equilibrium theory (see  \citealt{WANG2013equilibrium}).
The authors propose a metaheuristic algorithm where the multi-objective UL problem is addressed by NSGA-II and the LL problem by a quasi-Newton method, and solve an illustrative example, showing the conflicting nature of the three considered UL objectives.
In \cite{wen2016minimizing}, an illustrative example solved by simulating different toll pricing strategies and users' behavior reveals the impact of various constant and variable pricing strategies on travel time and CO\textsubscript{2} emissions.
The LL problem solves a dynamic user equilibrium, the simplest dynamic extension of Wardrop’s equilibrium, based on which all routes used by travelers departing at the same time with same origin and destination must have equal and minimal travel time (\citealt{SzetoWong}).
The authors rely on queuing theory for modeling the LL equilibrium.
In the problem proposed by \cite{zhao2016bi}, the leader aims at setting an optimal GHG emissions charge scheme for travelers (emission tolls), who decide on their travel mode (e.g., cars, motorcycles, bus). The solution method integrates an EA for the UL and the ‘‘all-or-nothing” (AON) algorithm (all the demand of each OD pair is assigned to the shortest path for that pair) for the LL traffic assignment. In the multimodal adaptation of LL equilibrium formulation \eqref{eq17eq}-\eqref{eq22eq}, there is an additional objective function term to consider the different travel cost for each travel mode times the related probability of choosing that mode.
The authors provide managerial insights from their numerical example, highlighting the trade-off between emission reduction goals and system efficiency, where travel time costs increase more rapidly with each additional percentage reduction in emissions.

\paragraph{Bike lanes design}
\label{par:bike}
Bike lanes are an alternative transport mode that indirectly reduce carbon emissions from vehicles. 
In bilevel bike lanes networks optimization problems there are multiple followers, the network users, interacting with the leader, the network designer, and with each other in a network equilibrium problem. 

For problems in which the equilibrium in the LL involves cyclists and motorists (i.e., multimodal user equilibrium traffic flow model), we mention \cite{sohn2011multi} and \cite{bagloee2016bicycle}.
\cite{sohn2011multi} deals with road diet network design problems (i.e., how to dedicate road portions exclusively to cyclists without affecting motorists).
He presents two multi-objective bilevel models, in which motorists can and cannot change their travel mode.
Cyclists travel times in the LL are assumed to be fixed once the road diet is decided in the UL. 
When travel modes are fixed, the motorists equilibrium in the LL is formulated as in model \eqref{eq17eq}-\eqref{eq22eq} optimizing automobile flow considering all arcs affected by the leader decision (road diet).
When travel modes are variable, the motorists equilibrium in the LL is adapted to include the mode choice.
A second term in the objective function \eqref{eq17eq} integrates the bicycle flows, ensuring that demand allocation (motorists/cyclists) respects the binary choice logit model and a new equilibrium is achieved, balancing flows across both modes while optimizing motorists' route choice.
NSGA-II algorithm is applied to a small test network to find an approximated Pareto front of non-dominated solutions for each model.
\cite{bagloee2016bicycle} tackle bike lanes network design in congested areas.
In the UL, binary design variables are used to select roads for exclusive bicycle use among the initial mixed-use road network.
In the LL, flow variables $y_a$ divide in the flows of motorists and cyclists ensuring multiclass equilibrium.
The LL is an extension of model \eqref{eq17eq}-\eqref{eq22eq}, where capacity and road diet coupling constraints are introduced such that the user equilibrium is constrained by the UL decision.
In this multiclass equilibrium, the feasibility and uniqueness of the solution are not guaranteed, as reported by the authors.
The authors present an effective B\&B algorithm, in which branching is iteratively applied to partial UL solutions until a complete one is found and multiclass user equilibrium is solved to find a new (incumbent) upper bound until no partial solutions are found. Their algorithm is able to solve a real-world instance with around 1000 nodes and 3000 arcs.
In \cite{gaspar2015bilevel}, the followers include car users, bus users, and cyclists and the LL problem is a modal-split assignment problem, in which users decide their optimal route and travel mode.
The proposed GA is applied to a real-world medium-size case study of the city of Santander (Spain) with 1500 candidate links for bicycle lanes.

In the bilevel bike lanes network design problem by \cite{LiuSJ22}, the leader designs the bike network to maximize total utility, while the followers (the cyclists) 
are in equilibrium with a multinomial logit model for route choices. The latter is a different formulation in which the probability $p_{mr}$ for cyclist $m$ of choosing route $r$ is the LL variable and the equilibrium model finds the minimum cost (maximum net utility) choice for every cyclist. The objective function is the sum of a convex term (balanced distribution of probabilities) and a linear term (total utility).
The authors first linearize the LL problem and then derive the KKTR, which is proven to be asymptotically exact. The authors are able to solve a real-world large-scale case study in China with Gurobi. They also show how the introduction of cycling route choices (the bilevel model) impacts the solution with respect to the single-level model, which only maximizes cyclist's utility function.

Bike lanes are also connected to sidewalks and crosswalks, which are designed similarly with bilevel optimization.
\cite{rashidi2016optimal} explore pedestrian safety in traffic management through bilevel optimization. They formulate an SLMF bilevel problem in which the city planner acts as the leader, making decisions on the location of sidewalks and crosswalks. The travelers, as followers, decide on their transportation mode (i.e., automobile, public transport, or walking) and routes. The LL problem is 
represented using travel cost equilibrium constraints in the UL.
The authors present a non-convex SLR, drawing from existing literature to certify the existence and uniqueness of the lower-level solution based on the nonlinear complementary conditions. However, in practice, the model struggles to generate feasible bilevel solutions in reasonable time.
Two heuristic methods calling the same nonlinear complementary algorithm to solve the LL equilibrium are presented to solve small illustrative examples, showing that pedestrian safety can be improved and costs can be reduced.

\paragraph{Green facility location}
\label{par:flp}
In the FLP, the objective is to strategically position a set of facilities to minimize the cost of meeting the demands of a group of customers (see \citealt{laporte2019facility}).
In the bilevel FLPs, the leader typically decides on the location and size of the facilities that could be used by the followers.
In the cases where the followers are drivers optimizing routes and facility choices, the LL problem is modeled as an equilibrium problem.
In green bilevel FLPs, environmental or social decisions and goals are tackled by the leaders.
\cite{ChenQMU20} formulate the bilevel problem of optimal location and capacity planning of electric vehicle (EV) charging stations. 
The classical LL equilibrium in model \eqref{eq17eq}-\eqref{eq22eq} is adapted here to consider the total flow on links as the sum of gasoline vehicles and EVs arrivals. In the equilibrium, the travel and waiting time at EV charging stations is minimized for every follower (charging and route choice equilibrium), also considering coupling constraints that restrict flow assignments to only constructed facilities.
The authors solve the proposed KKTR with an iterative algorithm (IA), which certifies the C-stationarity of bilevel solutions, on different illustrative examples (up to around 30 nodes, 80 links, and 40 candidate locations). 
\cite{zhou2023location} consider a bilevel model with a multimodal network equilibrium in the LL for the optimal placement of transit-oriented development (TOD) stations in a Chinese city. TOD is an urban planning approach focused on creating compact and mixed-use urban areas closely linked with mass transit stations to efficiently integrate jobs, housing, and services. 
The LL is a multimodal stochastic user equilibrium (SUE). The SUE differs from Wardrop equilibrium in the fact that followers make probabilistic mode choices (logit model) to minimize perceived travel times.
The authors develop a metaheuristic algorithm and solve a large-size real-world instance. They also show how TOD stations can reduce carbon emissions and how they impact multimodal networks.

In other cases, the leader is an authority (e.g., the government) and the followers are one or more individuals.
\cite{GangTLXSY15} include green goals in an FLP for stone industrial parks, which are districts of stone enterprises located close to each other to share infrastructures and control costs and pollution (from, e.g., dust and water consumption).
Several sustainable bilevel applications deal with industrial parks optimization problems. We dedicate below a specific section to these problems (see Section \ref{subsec:ind_eco}).
In the FLP studied by \cite{GangTLXSY15}, the leader in the UL is the local government taking parks location decisions ($x$), whereas the followers in the LL are the stone enterprises in conflict with each other but in cooperation with the leader to decide on park selection ($y$) subject to a constraint on maximum emission and coupling capacity constraints as in constraints \eqref{eq19} of model \eqref{eq16}-\eqref{eq20} for UL/LL network design variables. In this case, however, there is no pricing but pure network design in the UL with two linear objective functions ($F_1(y) = a_1y$ and $F_2(x,y) = a_2x + b_2y$).
The authors employ a matheuristic co-evolutionary approach in which UL and LL exchange information gradually in the iterative process. They develop a PSO algorithm based on the concept of satisfactory level of objectives for both the UL and LL and solve a case study in China, showing a 20\% reduction for all the objectives of both the leader and the followers.

\paragraph{Air passenger transportation}
\label{par:air}
\cite{qiu2020carbon} present an SLMF bilevel model for determining the optimal carbon tax incentive policy for air passenger transportation fossil fuel saving and carbon emissions reduction. 
The leader is the government, setting the carbon tax and allocating carbon tax incentive subsidies to the airlines for maximizing the net social benefit of the policy and to increase the system efficiency. Airlines are the followers, developing their air transport planning to maximize their profit.
Their mathematical model refers to price setting formulation \eqref{eq8}-\eqref{eq11}, with UL variable $T$ being the government carbon tax price to airlines, and LL variables $y$ and $z$ being the improvement (e.g., fuel type, fuel efficiency) and adjustment (tax transfer to passengers) strategies to adapt to UL price.
A GA is proposed to solve a case study in China. Results show that a carbon tax incentive policy has the potential to save fuel and improve environment under proper conditions.

\subsection{Sustainable routing}
\label{subsec:rout}
In the VRP, routes must be optimized to visit a set of customers with minimum total travel cost or time (see \citealt{tothvigo2014}).
Traditionally, these objectives have been purely economic. However, in recent years, the classical VRP has been widely extended to directly incorporate environmental and social impacts, moving beyond solely economic goals (\citealt{BEKTAS20111232}). 
This analysis focuses specifically on these sustainable extensions of the VRP.
While several studies have explored sustainability-driven variations of VRP, relatively few, to the best of our knowledge, have approached them using bilevel optimization (see Table \ref{tab:routing}).

The Pollution Routing Problem (PRP) considers customer time windows, vehicle speeds and loads, and aims to consider carbon emissions in the road freight transport sector (\citealt{BEKTAS20111232}).
The classical PRP objective functions are the minimization of total cost of fuel consumption and GHG emissions, total traveled distance, total weighted load, or total energy consumed by vehicles. 
Very few PRPs in the literature show a bilevel formulation.
\cite{nath2018novel} study a bilevel multi-objective PRP where the depot acts as the leader and vehicles as followers. In the UL, the customers are assigned to vehicles to minimize total traveled distance and fuel consumption, whereas, in the LL, vehicle routes are determined by minimizing distance.
The authors employ a metaheuristic algorithm based on NSGA-II for the UL and a GA for each vehicle in the LL, yielding improved results compared to NSGA-II applied to a single-level multi-objective PRP. \cite{qiu2020carbonpricing} present another bilevel PRP, incorporating a carbon pricing system.
The authors propose a price setting formulation (see \eqref{eq8}-\eqref{eq11}).
The authority acts as the leader and sets a carbon tax and a carbon subsidy to a freight transport company based on actual carbon emissions and emissions reductions, respectively.
The company acts as the follower making routing decisions. 
The LL appears as a classical PRP with total costs minimization, including the driver $(c_2)$ and fuel $(c_1)$ costs, the cost to pay for the carbon tax $(T_1)$ set by the leader for emissions $(y)$, and the carbon subsidy $(T_2)$, which is a positive cost (penalty price) if actual emissions $(y)$ exceed initial emissions (constant $\bar{y})$ and negative (subsidy for decreasing emissions) on the contrary.
The objective function is an adaptation of equation \eqref{eq9} given as 
\begin{equation*}
    \min_{y}  \quad (c_1 + T_1)y + T_2(y-\bar{y})+ c_2z.
\end{equation*}
Carbon emissions are obtained by converting fuel consumptions, a function of LL decisions on routing, flow, service time, and speed. 
The problem is solved by means of an interactive heuristic method based on PSO and adaptive large neighborhood search for the UL and LL, respectively. 
Results based on real-world UK instances reveal that carbon pricing initiatives can significantly reduce emissions with a relatively modest increase in costs.

On the other hand, VRPs with social implications address public transportation systems, an area that is gaining increasing attention in VRP literature (see, e.g., \citealt{DELLEDONNE2025}), though still less extensively explored compared to economic and environmental aspects and very rarely approached with bilevel optimization.
Efficient public transport systems remain crucial for urban development, accessibility, sustainability, and overall societal well-being (\citealt{EUUrbanMobility2021}).
In this analysis, we highlight one notable bilevel application in the literature with direct social implications.
\cite{PARVASI2017} and \cite{calvetegale2023} study two closely related bilevel problems related to the school bus routing optimization with students' choice consideration. 
In \cite{PARVASI2017}, the leader in the UL is the public transportation company deciding bus stations location and bus routes. The students are the independent followers in the LL problems, who decide whether to take which bus and at which station. The authors propose two hybrid metaheuristic algorithms based on GA, simulated annealing (SA), and tabu search (TS) effectively applied on random instances.
In \cite{calvetegale2023}, the leader of the SLMF bilevel problem is the authority who decides the routes, given the bus stops, and the number of students assigned to each stop considering students/stops accessibility and preferences.
The students, in the LL problem, decide which bus stop to select. The authors provide an SLR based on duality theory, which works on small-size instances with optimality certificate, and a metaheuristic algorithm, which obtains good performances on large-size benchmark instances.

\afterpage{
\clearpage
\newgeometry{margin=3cm}
\begin{landscape}
\centering
\begin{table}[ht]
    \tiny
    \centering
    \caption{Summary of papers on bilevel applications focused on carbon policies on production planning. 
    }
    \label{tab:production-pol}
    \begin{tabularx}{\linewidth}{p{0.10\textwidth}p{0.17\textwidth}p{0.05\textwidth}p{0.02\textwidth}p{0.15\textwidth}p{0.15\textwidth}p{0.22\textwidth}p{0.22\textwidth}X}
    \toprule
    Application & References & \multicolumn{2}{c}{Bilevel game} & \multicolumn{2}{c}{Players} & \multicolumn{2}{c}{Objectives} & Method\\
    \cmidrule(lr){3-4}\cmidrule(lr){5-6}\cmidrule(lr){7-8}
    & & Class & Eq. & UL & LL & UL & LL &  \\
    \cmidrule(lr){1-1}\cmidrule(lr){2-2}\cmidrule(lr){3-4}\cmidrule(lr){5-6}\cmidrule(lr){7-8} \cmidrule(lr){9-9}       
    \multirow{5}{*}{\parbox{0.10\textwidth}{Carbon policies}} 
    & \cite{sinha2013} & SLSF & & Government & Mining company & Max revenues; Min pollution & Max profit & Hybrid MOEA \\
    & \cite{almutairi2014carbon} & SLSF & & Government & Industry & Min production and target deviation weighted difference & Max surplus & KKTR \\
    & \cite{HONG2017172} & SLMF & \checkmark & Local government & Firms & Max social welfare & Max profit & DP + BSA + GA \\  
    & \cite{Molavi2020} & SLMF & & Regulatory authority & Ports & Min emission caps & Max profit & HEUR \\
    & \cite{singh2018} & SLSF & & Manufacturer/Seller & Seller/Manufacturer & Max profit & Max profit & KKTR \\
    \bottomrule
    \end{tabularx}
\end{table}
\vspace{-.6cm}
    \begin{table}[H]
    \tiny
    \centering
    \caption{Summary of papers on bilevel applications focused on green product design and manufacturing.  
    }
    \label{tab:production-design}
    \begin{tabularx}{\linewidth}{p{0.10\textwidth}p{0.17\textwidth}p{0.05\textwidth}p{0.02\textwidth}p{0.15\textwidth}p{0.15\textwidth}p{0.22\textwidth}p{0.22\textwidth}X}
    \toprule
    Application & References & \multicolumn{2}{c}{Bilevel game} & \multicolumn{2}{c}{Players} & \multicolumn{2}{c}{Objectives} & Method\\
    \cmidrule(lr){3-4}\cmidrule(lr){5-6}\cmidrule(lr){7-8}
    & & Class & Eq. & UL & LL & UL & LL &  \\
    \cmidrule(lr){1-1}\cmidrule(lr){2-2}\cmidrule(lr){3-4}\cmidrule(lr){5-6}\cmidrule(lr){7-8} \cmidrule(lr){9-9}       
    \multirow{4}{*}{\parbox{0.10\textwidth}{Green product design}} 
    & \cite{machencai2018} & SLSF & & Manufacturer & Supplier & Max utility/cost ratio & Min total costs & GA \\
    & \cite{ZHU2020114134} & SLSF & & Designer & Designer & Min fuel consump.; Min GHG em.; Min net present cost & Min fuel consump. & MOPSO + HEUR \\
    &  \cite{wujia2020} & SLMF & & Platform operator & Service demanders & Max environmental, economic, social perf. &  Min time; Min price; Max quality & Hybrid MOPSO \\
    & \cite{zengdong2020} & SLMF & & Station owner & Station users & Max net revenue & Max utility & KKTR + CCG \\
    \bottomrule
    \end{tabularx}
\end{table}
\vspace{-.6cm}
\begin{table}[H]
    \tiny
    \centering
    \caption{Summary of papers on bilevel applications focused on water management. 
    }
    \label{tab:water-man}
    \begin{tabularx}{\linewidth}{p{0.10\textwidth}p{0.17\textwidth}p{0.05\textwidth}p{0.02\textwidth}p{0.15\textwidth}p{0.15\textwidth}p{0.22\textwidth}p{0.22\textwidth}X}
    \toprule
    Application & References & \multicolumn{2}{c}{Bilevel game} & \multicolumn{2}{c}{Players} & \multicolumn{2}{c}{Objectives} & Method\\
    \cmidrule(lr){3-4}\cmidrule(lr){5-6}\cmidrule(lr){7-8}
    & & Class & Eq. & UL & LL & UL & LL &  \\
    \cmidrule(lr){1-1}\cmidrule(lr){2-2}\cmidrule(lr){3-4}\cmidrule(lr){5-6}\cmidrule(lr){7-8} \cmidrule(lr){9-9}       
    \multirow{2}{*}{\parbox{0.10\textwidth}{Water scarcity}} & \cite{ANANDALINGAM1991233} & SLMF  &  & United Nations & Countries & Max economic net benefit & Max economic net benefit & Dual. SLR \\
    & \cite{calvete2023} & SLMF & & Government & Water demand points managers & Min water deficit, Max satisfaction, Min water price & Max economic net benefit & KKTR+Lexic. \\
    \cmidrule(lr){1-1}\cmidrule(lr){2-2}\cmidrule(lr){3-4}\cmidrule(lr){5-6}\cmidrule(lr){7-8} \cmidrule(lr){9-9}
    \multirow{1}{*}{\parbox{0.10\textwidth}{Water pollution}} & \cite{ZHAO2013377} & SLMF  &  & Lake authority & Local authorities & Min pollution reduction cost & Min pollution and transfer cost & KKTR \\
    \bottomrule
    \end{tabularx}
\end{table}
\vspace{-.6cm}
\begin{table}[H]
    \tiny
    \centering
    \caption{Summary of papers on bilevel applications focused on waste management. 
    }
    \label{tab:waste-man}
    \begin{tabularx}{\linewidth}{p{0.10\textwidth}p{0.17\textwidth}p{0.05\textwidth}p{0.02\textwidth}p{0.15\textwidth}p{0.15\textwidth}p{0.22\textwidth}p{0.22\textwidth}X}
    \toprule
    Application & References & \multicolumn{2}{c}{Bilevel game} & \multicolumn{2}{c}{Players} & \multicolumn{2}{c}{Objectives} & Method\\
    \cmidrule(lr){3-4}\cmidrule(lr){5-6}\cmidrule(lr){7-8}
    & & Class & Eq. & UL & LL & UL & LL &  \\
    \cmidrule(lr){1-1}\cmidrule(lr){2-2}\cmidrule(lr){3-4}\cmidrule(lr){5-6}\cmidrule(lr){7-8} \cmidrule(lr){9-9}       
    Urban waste & \cite{CARAMIA2022} & SLSF  &  & Central authority & Local authority & Min land-use stress, Min impact on health & Min transportation cost & KKTR+WS \\
    \cmidrule(lr){1-1}\cmidrule(lr){2-2}\cmidrule(lr){3-4}\cmidrule(lr){5-6}\cmidrule(lr){7-8} \cmidrule(lr){9-9}
    Hazardous waste & \cite{AMOUZEGAR1999100} & SLMF  &  & Central authority & Firms & Min total costs & Min cost & HEUR \\
    \bottomrule
    \end{tabularx}
\end{table}
\end{landscape}
\restoregeometry
}

\section{Production planning and manufacturing}
\label{sec:production}
Governments' environmental policies, such as carbon taxes, emissions regulations, and renewable energy incentives, along with public awareness among stakeholders like manufacturers, energy producers, and transportation companies, are driving the adoption of green technologies and sustainable industrial practices. Stakeholders are encouraged or required to adopt renewable energy, energy-efficient technologies, and advanced pollution control systems to meet regulations and public expectations
(see Tables \ref{tab:production-pol} and \ref{tab:production-design}). 
These efforts align with SDG 8 (Decent Work and Economic Growth) by fostering a green economy, SDG 9 (Industry, Innovation, and Infrastructure) through manufacturing growth, economic acceleration, and CO\textsubscript{2} emissions reduction, and SDG 12 (Responsible Consumption and Production) by promoting policies for sustainable consumption of stakeholders and carbon footprint reduction of industrial activities.
In bilevel frameworks, either the government imposes a carbon policy and the industry responds or multiple industrial entities interact with each other to reduce the environmental impact of industrial outputs, typically in the UL, while seeking efficiency in both UL and LL.

\subsection{Carbon policies}
\label{subsec:carbonpol}
Governmental policies and industrial reactions can be generally modeled as bilevel nonlinear price setting problems. 
Carbon policies to control industrial activities pollution are mainly direct and indirect taxes, subsidies, and tradable permits. 
In most bilevel applications, the leader in the UL is the government imposing carbon taxes and one or more firms are the followers taking operational decisions in the LL.
\cite{sinha2013} study an SLSF bilevel problem to optimize government's environmental regulatory decisions on a mining project. The UL decision is on the optimal tax rate for each time period to be paid by the mining company. The LL decision of the company is on the extraction amount for each time period. The authors propose a hybrid multi-objective evolutionary algorithm in a nested framework, where for each UL decision the follower's decision is approximated under an optimistic approach, and solve a case study in Finland, providing a valuable approximate revenues/pollution Pareto front of solutions.
In a similar SLSF bilevel model studied by \cite{almutairi2014carbon}, 
the follower is a generic industry, determining production quantities and selecting fuel types.
The authors present a KKTR approach and apply it to solve a cement industry case study in Canada to optimality. Results show that moderate tax rates can bring large reductions in the overall emissions and the optimal tax rate reduces the quantity supplied by high emitters and increases the quantity supplied by low and average emitters.
For the SLMF setting of multiple companies, \cite{HONG2017172} introduce a bilevel model for a carbon emissions cap-and-trade scheme. The local government sets the emission trading scheme and 
firms make production planning and technology selection decisions and trading emission allowances among each other.
The followers' problem is modeled as a Cournot equilibrium model. The authors propose a hybrid algorithm that combines dynamic programming (DP) for followers' decision, a binary search algorithm (BSA) for the LL equilibrium to obtain an $\epsilon$-approximate optimal solution, and a GA for the leader's decision. The effectiveness of the cap-and-trade scheme is shown through a numerical example in China.
\cite{Molavi2020} model an SLMF problem where a regulatory authority, as leader, minimizes a weighted sum of emission caps, while ports, as independent followers, maximize profits considering investment costs, carbon taxes, and carbon incentives. Using a tailored heuristic based on upper and lower bounds, they solve a case study on US ports, finding that taxes alone reduce emissions, but combining taxes with incentives enhances reductions without harming port profitability.
Recent studies use fuzzy optimization to approximate bilevel solutions in government-industry carbon policy problems (see, e,g., \citealt{TAN2023}). 

In other contexts like the one studied by \cite{singh2018}, the bilevel game is between two firms and the government sets green taxes and duties from outside. The authors present two SLSF bilevel models for off-shore manufacturing contracts optimization with pollution control. In the first model, the leader is a seller firm in a developed country deciding the retail price of a product in its market, and the follower is another firm in a developing country that manufactures the product with a lower manufacturing cost and decides on the transfer price to the seller. 
A green tax and an import duty are paid, respectively, by the manufacturer and the seller to their governments. In the second model, on the contrary, the leader is the manufacturer and the follower is the seller. The authors obtain KKTRs, linearized and solved to optimality. A case study on electronics between a seller in USA and a manufacturer in China shows that the bilevel models allow to improve the contract for both parties under pollution control policies.

\subsection{Green product design and manufacturing}
\label{subsec:greentechn}
Green product design and manufacturing practices are spreading across different industrial sectors and various bilevel models are consequently emerging. 
In SLSF problems, two entities of the same supply chain interact with each other with sustainable goals (see Table \ref{tab:production-design}).
In \cite{machencai2018}, 
the manufacturer is the leader, deciding on production modules, paying the carbon tax, and giving carbon incentives to reduce emissions. The supplier acts as the follower, deciding on product modules production.
The authors propose a nested GA and show an illustrative example on notebook computers, showing the effectiveness of their carbon incentive-based bilevel model in improving both the product profit and the carbon emission.
In \cite{ZHU2020114134},
the leader is the designer of hybrid EVs 
and the follower is another EV designer making decisions on the energy management system.
The proposed solution method combines MOPSO at the UL and a heuristic rule-based control approach at the LL. The results on an illustrative example show that the bilevel approach is able to find optimal solutions resulting in less fuel consumption, less GHG emissions, and less net present cost compared to the single-level optimization.

In SLMF games, the followers are the users of a green product or service.
\cite{wujia2020} propose a bilevel model to assess the three dimensions of sustainability for cloud manufacturing services:  environmental performance (e.g., energy consumption), economic performance (e.g., flexibility), and  social performance (e.g., consumer satisfaction).
In the UL, the platform operator decides which demanded services to fulfill and how.
In the LL, the service demanders determine the operational strategy for service composition.
The authors propose a hybrid PSO algorithm and show that it is more efficient than the classical PSO on the single-level multi-objective model. 
\cite{zengdong2020} introduce a robust SLMF bilevel model for the design and planning of a public Plug-in EVs charging station.
The leader, the station owner, decides on the renewable energy and energy storage capacity and retail price.
The followers, the station users, decide on their charging-energy demand, after uncertainty on renewable availability, wholesale energy price, and number of EVs arrivals realizes. The authors propose a KKTR for the model and a column-and-constraint generation (CCG) algorithm to solve a case study. They compare the results of the bilevel approach with a single-level one and show the benefits of the former to increase renewable availability and decrease energy costs. 

\section{Water, waste, and agriculture management}
\label{sec:waste}
We explore a selection of works leveraging bilevel optimization in water, waste, and agriculture management.
These studies align with SDG 3 (Good Health and Well-being) by mitigating pollution-related health risks.
They contribute to SDG 6 (Clean Water and Sanitation) by addressing water scarcity and pollution, ensuring sustainable resource use, and emphasizing the need to combat lake ecological deterioration, a growing global concern. SDG 15 (Life on Land) is supported through sustainable agriculture, focusing on preventing land degradation and minimizing the environmental impact.

\subsection{Water management}
\label{subsec:water}
Water is a primary source of life and its management has consistently posed complex environmental and societal challenges. Multiple stakeholders are involved in the water management process and different issues may arise. We cover water scarcity and water pollution (see Table \ref{tab:water-man}). 
\afterpage{
\clearpage
\newgeometry{margin=3cm}
\begin{landscape}
\begin{table}[ht]
    \tiny
    \centering
    \caption{Summary of papers on bilevel applications focused on agriculture. 
    }
    \label{tab:agri}
    \begin{tabularx}{\linewidth}{p{0.10\textwidth}p{0.17\textwidth}p{0.05\textwidth}p{0.02\textwidth}p{0.15\textwidth}p{0.15\textwidth}p{0.22\textwidth}p{0.22\textwidth}X}
    \toprule
    Application & References & \multicolumn{2}{c}{Bilevel game} & \multicolumn{2}{c}{Players} & \multicolumn{2}{c}{Objectives} & Method\\
    \cmidrule(lr){3-4}\cmidrule(lr){5-6}\cmidrule(lr){7-8}
    & & Class & Eq. & UL & LL & UL & LL &  \\
    \cmidrule(lr){1-1}\cmidrule(lr){2-2}\cmidrule(lr){3-4}\cmidrule(lr){5-6}\cmidrule(lr){7-8} \cmidrule(lr){9-9}       
    \multirow{4}{*}{\parbox{0.10\textwidth}{Agricultural policy instrument}} & \cite{candlernorton1977} & SLMF & & Policy maker & Producers & Max employment, income, production & Max profit & Exact \\
    & \cite{ONAL1995227} & SLMF & \checkmark & Policy maker & Farmers & Max total revenue & Max surplus & KKTR+HEUR \\
    & \cite{WHITTAKER201715} & SLMF & & Policy maker & Farmers & Max total profit, Min nitrogen in land & Max profit & Hybrid GA \\
    & \cite{barnhart2017} & SLMF & & Policy maker & Farmers & Max total profit, Min fertilizer pollution & Max profit & MOEA \\
    \cmidrule(lr){1-1}\cmidrule(lr){2-2}\cmidrule(lr){3-4}\cmidrule(lr){5-6}\cmidrule(lr){7-8} \cmidrule(lr){9-9}  
    \multirow{1}{*}{\parbox{0.10\textwidth}{Agricultural supply chain}} & \cite{albornoz2023} & SLSF & & Producer & Wholesaler & Max profit & Min cost & KKTR \\
    \bottomrule
    \end{tabularx}
\end{table}
\vspace{-.6cm}
    \begin{table}[H]
    \tiny
    \centering
    \caption{Summary of papers on bilevel applications focused on green supply chain design and planning. 
    }
    \label{tab:sc}
    \begin{tabularx}{\linewidth}{p{0.10\textwidth}p{0.17\textwidth}p{0.05\textwidth}p{0.02\textwidth}p{0.15\textwidth}p{0.15\textwidth}p{0.22\textwidth}p{0.22\textwidth}X}
    \toprule
    Application & References & \multicolumn{2}{c}{Bilevel game} & \multicolumn{2}{c}{Players} & \multicolumn{2}{c}{Objectives} & Method\\
    \cmidrule(lr){3-4}\cmidrule(lr){5-6}\cmidrule(lr){7-8}
    & & Class & Eq. & UL & LL & UL & LL &  \\
    \cmidrule(lr){1-1}\cmidrule(lr){2-2}\cmidrule(lr){3-4}\cmidrule(lr){5-6}\cmidrule(lr){7-8} \cmidrule(lr){9-9}       
    \multirow{7}{*}{\parbox{0.10\textwidth}{Green supply chain}} & \cite{avval2023joint} & SLMF & & Government & SCs decision makers & Min carbon cap & Min cost & ILS \\
    & \cite{camacho2022tabu} & SLSF & & SC distributor & SC manufacturer & Max profit; Min CO\textsubscript{2} emissions & Min cost & HEUR \\
    & \cite{camacho2023hierarchized} & SLSF & & SC distributor & SC manufacturer & Max profit; Min CO\textsubscript{2} emissions & Min cost & HEUR \\
    & \cite{cantu2021novel} & SLSF & & SC designer & SC operational decision maker & Min cost; Min GWP & Min cost; Min GWP & MOEA+LP \\
    & \cite{cantu2023capturing}  & SLSF & & SC designer & SC operational decision maker & Min cost; Min GWP & Min cost; Min GWP & MOEA+LP \\
    & \cite{ghomi2021competitive} & SLSF & & SC1 & SC2 & Max profit cost; Min CO\textsubscript{2} emissions &  Max profit & SLR+$\epsilon$-constr. \\
    & \cite{golpira2017robust} & SLSF & & SC1 & SC2 & Max profit cost; Min CO\textsubscript{2} emissions &  Max profit & SLR+$\epsilon$-constr. \\
    \bottomrule
    \end{tabularx}
\end{table}
\vspace{-.6cm}
\begin{table}[H]
    \tiny
    \centering
    \caption{Summary of papers on bilevel applications focused on Eco-Industrial Parks (EIPs). 
    }
    \label{tab:eip}
    \begin{tabularx}{\linewidth}{p{0.10\textwidth}p{0.17\textwidth}p{0.05\textwidth}p{0.02\textwidth}p{0.15\textwidth}p{0.15\textwidth}p{0.22\textwidth}p{0.22\textwidth}X}
    \toprule
    Application & References & \multicolumn{2}{c}{Bilevel game} & \multicolumn{2}{c}{Players} & \multicolumn{2}{c}{Objectives} & Method\\
    \cmidrule(lr){3-4}\cmidrule(lr){5-6}\cmidrule(lr){7-8}
    & & Class & Eq. & UL & LL & UL & LL &  \\
    \cmidrule(lr){1-1}\cmidrule(lr){2-2}\cmidrule(lr){3-4}\cmidrule(lr){5-6}\cmidrule(lr){7-8} \cmidrule(lr){9-9}       
    \multirow{2}{*}{\parbox{0.10\textwidth}{Materials, water, and energy exchange}} & \cite{ramos2016water};\cite{ramos2016optimal} & SLMF & & EIP authority & EIP companies & Min total water consumption & Min cost & KKTR \\
    & & MLSF & & EIP participants & EIP authority & Min cost & Min total water consumption & KKTR \\
    \cmidrule(lr){1-1}\cmidrule(lr){2-2}\cmidrule(lr){3-4}\cmidrule(lr){5-6}\cmidrule(lr){7-8} \cmidrule(lr){9-9}  
    \multirow{2}{*}{\parbox{0.10\textwidth}{Utility share}} & \cite{ramos2018utility} & SLMF & & EIP authority & EIP participants & Min total CO\textsubscript{2} emissions & Min cost & KKTR \\
    & & MLSF & & EIP participants & EIP authority & Min cost & Min total CO\textsubscript{2} emissions & KKTR \\
    \cmidrule(lr){1-1}\cmidrule(lr){2-2}\cmidrule(lr){3-4}\cmidrule(lr){5-6}\cmidrule(lr){7-8} \cmidrule(lr){9-9}  
    \multirow{1}{*}{\parbox{0.10\textwidth}{Carbon incentives}} & \cite{gu2020bi} & SLMF & & EIP authority &  EIP participants & Max revenue & Min cost & Iterative primal dual \\
    \bottomrule
    \end{tabularx}
\end{table}
\end{landscape}
\restoregeometry
}

\paragraph{Water scarcity}
\label{par:water-shortage}
Bilevel optimization has been applied for water conflicts resolution.
\cite{ANANDALINGAM1991233} introduce an insightful problem in international rivers. They present a comparative analysis of different hierarchies in multi-level Stackelberg games involving an arbitrator (specifically, the United Nations) and two countries (India and Bangladesh) who compete for the use of Ganges river waters for hydroelectric power, irrigation, and flood protection. 
In the bilevel hierarchy, the arbitrator is the leader and the two countries act as followers of equal status, in conflict with both each other and the leader. 
In the three-level hierarchy, the arbitrator again assumes the leader role, but the two followers are engaged in a leader-follower hierarchical relationship. 
Both models are solved exactly by SLRs and using the penalty function method, given strong duality conditions on the LL problem.
The authors show that higher objective function values are obtained under the three-level model, owing to substantial subsidies from the arbitrator covering a significant portion of the system costs.
\cite{calvete2023} present an SLMF bilevel formulation to address the optimal allocation of water to demand points that are in conflict with each other due to water shortage. The leader is the government deciding on water allocation to demand points and establishing fees. In the LL, the multiple followers are demand points managers who make decisions regarding water distribution within their areas. 
The bilevel model is reformulated as a single lexicographic multi-objective MILP model by the KKTR and tested effectively on several water system scenarios.

\paragraph{Water pollution}
\label{par:water-pollution}
In the work of \cite{ZHAO2013377}, an SLMF bilevel programming formulation is employed to model the hierarchical structure of authorities responsible for lake water pollution control. In the UL, there is one leader, the lake authority. In the LL, there are multiple followers, the sub-regional authorities. 
The leader strategically sets tax rates while ensuring compliance with environmental quality standards. 
If the pollutant amount generated by a region exceeds the national standard, the region pays a fee proportional to the pollutant quantity.
Pollution reduction costs are set by the lake authority.
The followers, in turn, make decisions regarding their pollution reduction tactics.
The authors derive the KKTR and then apply the bisection method to find an approximated solution to the bilevel problem, providing a global convergence rate.

\subsection{Waste management}
\label{subsec:waste}
Authorities and citizens are compelled to assess the economic efficiency and environmental impact of waste due to its global volume increase and the high pollution rate of its treatment process.  
Within the waste management process, the national or regional authority is responsible for strategic decisions such as facility location for waste treatment, and local authorities or private companies deal with operational decisions such as routing for waste collection (see Table \ref{tab:waste-man}).

\paragraph{Urban waste}
\label{par:urban-waste}
Municipal solid waste management encompasses a set of practices aimed at collecting and recycling solid waste generated by citizens, involving a complex network of facilities and related stakeholders, which must be efficiently designed.
In \cite{CARAMIA2022}, the global and local authorities serve as the leader and the follower responsible for waste network design and transportation decisions, respectively.
The problem involves location, capacity, and transportation decisions for collection centers, sorting facilities, landfills, and incinerators. 
The authors present two bilevel formulations under both the optimistic and pessimistic assumptions, along with corresponding KKTRs with the weighted sum (WS) method for dealing with the two UL objectives. They apply the linearized formulations to a case study in Thailand, showing that the results of bilevel programming are more realistic than the ones obtained with a bi-objective single-level approach.

\paragraph{Hazardous waste}
\label{par:haz-waste}
Hazardous waste includes all kinds of dangerous and ``special'' sources of waste, typically generated by industrial entities, such as chemicals and electronics.
Motivated by a case study in California, \cite{AMOUZEGAR1999100} introduce a capacity allocation and facility location SLMF bilevel model for hazardous waste management.
In their model, the central authority is the leader, who controls allocation and location decisions by setting prices and taxes for harmful policies, while the firms engaged in the waste management process act as followers.
The authors propose a two-phase heuristic algorithm that initially solves a relaxed single-level model, excluding the followers' objective, and subsequently solves the followers' problem. This penalty algorithm was proven to provide bilevel feasible solutions in \cite{Amouzegar1998}.

\subsection{Agriculture}
\label{subsec:agri}
Bilevel programming finds several applications in agriculture, particularly in public policy-making and environmental policies. Several studies focus on optimizing policy instruments, such as taxes, credits, and incentives, with environmental and social goals, anticipating the responses of farmers.
The bilevel problems are typically SLMF  and are solved by means of SLRs or heuristic algorithms (see Table \ref{tab:agri}).

\paragraph{Agricultural policy instruments}
\label{par:agr-poli}
In their pioneering work, \cite{candlernorton1977} define multi-level programming models in the economic policy context. They call ``policy problem'' and ``behavioral problem'', respectively, the UL and LL problems. The former represents the decision problem of policy makers who set policy instruments (e.g., tax rates, incentives). The latter is the decision problem of decentralized agents (e.g., firms, consumers, households) optimizing their economic behavior.
The authors propose a linear bilevel model application for a Mexican agricultural production region, where policy makers set subsidies, prices, budgets and taxes on materials and water, whereas the decentralized agents decide on production.
The application (solved by a ``hand-made'' algorithm inspired by the simplex method) illustrates the trade-offs between employment and production.

Other applications in agriculture model bilevel games between a policy maker and farmers.
\cite{ONAL1995227} address the optimal allocation of agricultural credits to farm groups in Indonesia. The leader in the UL allocates credits in a selected area. The LL problem is a market equilibrium problem for the farmers.
The KKTR of the bilevel problem is solved using a heuristic method based on the penalty function with bilevel feasibility certification. Results on a small real-world application show how reallocating credits, in particular decreasing credits allocation to large farmers in favor of small farmers, can pursue both agricultural market growth and equity.
Similar conclusions are shown by \cite{WHITTAKER201715}, who study the problem of setting green tax levels on fertilizer use in lands and show that bilevel optimization is effective for geographical targeting agri-environmental policies. 
The policy maker determines tax rates, whereas farmers make production decisions based on tax rates. The study proposes a bilevel optimization approach solved by a hybrid GA and applies it to the case study of Calapooia watershed (Oregon). Results open up to the possibility of including additional social objectives in the multi-objective UL.
\cite{barnhart2017} target a problem for the Raccoon watershed (Iowa) and study a deterministic problem with the assumption of the same tax rate for all farmers. They propose three EAs to identify the Pareto-optimal set of policies for the policy maker.
Their metaheuristics are based on the exact (for NSGA-II) or approximate optimal response of the followers, written as a closed-form solution in terms of UL decision variables.
The best method provides a well-established Pareto front for a large real-world instance with more than 1000 farmers. The authors also study the robust version of the problem considering production uncertainty, and produce more realistic Pareto fronts at the expense of the leader's objectives. 

\paragraph{Agricultural supply chain optimization}
\label{par:agr-sc}
\cite{albornoz2023} study a stochastic SLSF bilevel problem to address the selective harvest planning problem in a supply chain. 
The producer serves as the leader, making decisions related to harvest planning and scheduling (including harvest zone selection, quantities, and workforce planning). The wholesaler acts as the follower and determines the quantities to purchase from the producer.
Stochastic yields, demand, and prices are considered in the model. The KKTR of the stochastic model is used to solve a real-world instance for grape harvesting in Chile. Results show that the bilevel model leads to a 10\% increase in producer profits compared to a two-stage method where two separate models are solved in a hierarchical perspective. 

\section{Supply chains}
\label{sec:sc}
Supply Chain Network Design (SCND) is the strategic field related to the comprehensive planning and structuring of a supply chain.
In recent years, the field has evolved to incorporate environmental and social concerns, and the literature has studied several bilevel applications for sustainable SCND.
Studies align with SDG 7 (Affordable and Clean Energy) by promoting sustainable and reliable energy use in supply chains, SDG 8 (Decent Work and Economic Growth) by fostering the green economy and technological development, as well as SDG 9 (Industry, Innovation, and Infrastructure) and SDG 13 (Climate Action) by advancing resilient infrastructures, industrial growth, and CO\textsubscript{2} emissions reduction. Additionally, SDG 11 (Sustainable Cities and Communities) is addressed through reduced air pollution in supply chains and industrial parks planning.
We review applications of green SC design and planning (see Table \ref{tab:sc}) and eco-industrial parks (see Table \ref{tab:eip}).

\subsection{Green supply chain design and planning}
\label{subsec:sc_des}
In certain scenarios, the government acts as the leader of the bilevel problem influencing the regulations of multiple SCs, whereas the followers are the SC decision makers, namely the managers of the firms that have production and inventory facilities and deliver manufactured products.
In \cite{avval2023joint}, an SLMF bilevel problem is used to establish a cap-and-trade system for SCs.
The government is the leader, setting a carbon cap and carbon allowances to encourage SCs to adopt green technologies and reduce carbon emissions.
The SCs decision makers act as followers in a Stackelberg game, taking decisions according to their carbon allowances and trading allowances among each other. 
An iterated local search (ILS) is employed to solve small instances of the problem.
The UL solution quality and computational efficiency of their method is compared to an exact approach in which the bisection method applied to UL variables in the LL is combined with an optimization solver to produce an $\epsilon$-approximated LL optimal solution. 

In other scenarios, different stakeholders within the same SC can be involved in a bilevel problem, with one acting as the leader and the other as the follower. 
In \cite{camacho2022tabu}, the leader is a distribution company and the follower is a manufacturing company of the same SC. 
The distribution company acquires commodities from the manufacturing company and makes decisions on customer selection, routes, and vehicle types, while the manufacturing company defines the production plan to not exceed a maximum pollution rate.
An optimistic approach is assumed and a bi-objective nested TS algorithm, in which the follower’s problem is optimally solved for each leader’s solution, is used to find well-approximated emissions/profit Pareto fronts for instances with up to 1000 customers, 7 plants, and 80 vehicles.
An analogous bilevel model with profit maximization and CO\textsubscript{2} emissions minimization goals in the UL is presented in \cite{camacho2023hierarchized} in the context of SC restructuring (i.e., manufacturing and service level reduction) after the COVID-19 pandemic and consequent business closures.
They propose a multi-start heuristic algorithm with a well-posed LL problem under a tailored optimistic approach (i.e., with respect to the environmental UL objective) and solve medium-size instances with up to 700 customers, 60 vehicles and a few manufacturing plants. Results emphasize the strategic significance of achieving a balance between the two UL objectives for the supply chain leader.

\cite{cantu2021novel} study a sustainable SCND problem involving energy source, production, transportation and storage decisions of a hydrogen SC. A bi-objective MILP model minimizing the total daily cost and the GWP of the SC is reformulated as an SLSF bilevel problem. In the UL, the SC designer makes facility location decisions, while production and transportation decisions are addressed in the LL by the operational decision maker of the SC. 
A hybrid nested solution strategy, combining an MOEA for the UL problem and an LP-based heuristic approximating the transportation problem and introducing a single WS objective function with random weights for the LL problem, is proposed and compared with the classical $\epsilon$-constraint method applied to the single-level bi-objective MILP model. 
The results on six real-world instances in France show that the bilevel method provides better approximated cost/emissions Pareto fronts of solutions. 
Two years later \cite{cantu2023capturing} generalize the problem by including additional characteristics and nonlinear investment cost functions.
An adaptation of their previous hybrid solution strategy is shown to be valid since all the nonlinearities are included in the UL problem. The robustness of the solution method is proved by tests on the same case study in France.

In the third bilevel SCND scenario we consider, two supply chains can be involved in a hierarchical relationship: one SC acts as the leader (SC1) and the other as the follower (SC2), with uncertainty typically considered in the problem formulation.
\cite{ghomi2021competitive} present a closed-loop SCND problem (i.e., SC with a reverse flow of products for recycle) with two SCs competing on the same product market in a bilevel framework under demand uncertainty and disruption. The leader decides on the facility locations, supplier selection, demand satisfaction, retail price, inventory management, and reverse material flow.
The follower decides on their demand satisfaction and retail price.
A bi-objective KKTR is given for the stochastic bilevel model, and, to deal with
uncertainty, it is enriched with an RO approach. The resulting model is solved with the $\epsilon$-constraint method and tested on a real-world case study in Tehran (Iran).
A similar setting and solution approach are also proposed in \cite{golpira2017robust} for a green opportunistic SCND problem.

Other SCND applications include relief logistics and are addressed in Section \ref{subsec-emergency}.

\subsection{Eco-industrial parks}
\label{subsec:ind_eco}
An ``Eco-Industrial Park'' (EIP) is defined as a system of companies that are located close enough to exchange materials and energy and that aim to build sustainable economic, ecological, and social relationships, like an SC of strictly connected members.
Since 2000, and more extensively since 2010, several studies have addressed EIPs optimization problems mainly by means of multi-objective optimization and game theory (see \citealt{boix2015optimization}). 
More recently, bilevel programming has been introduced in the field for modeling two primary types of cooperation in an EIP: (i) the exchange of materials, water, and energy and (ii) the sharing of water and energy units.

\cite{ramos2016water} introduce a bilevel optimization model for the design and optimization of industrial water networks. They propose two original multi-leader-follower games involving the EIP authority regulating water exchanges and the companies. 
In the SLMF game, the leader is the authority, while the followers are the companies.
In the MLSF game, the roles are reversed. In the former case, the priority is given to environmental concerns, while the latter prioritizes the economic benefit of the companies. 
The two bilevel problems are both tackled by solving their KKTRs using non-convex and nonlinear programming solvers.
Results on a small example with three companies show that the bilevel approach (in both SLMF and MLSF cases) is able to advantage all three companies, while the multi-objective optimization approach only favours two of them.
Similar results have been shown in \cite{ramos2016optimal}. 

For EIP utility share, we only mention the relevant work by \cite{ramos2018utility}, who propose multi-leader-follower games for the design of a heat exchange utility network and introduce the concept of environmental authority (the EIP authority with the goal of minimizing the equivalent CO\textsubscript{2} emission from utility consumption).
Again, the authors formulate the SLMF and MLSF games and then solve their KKTRs.

Incentives for emissions reduction and energy-saving practices are also gaining attraction within EIPs.
\cite{gu2020bi} introduce energy price incentives for economic and environmental benefits of the partners of an EIP under the Chinese real-time multi-energy system. An SLMF bilevel model is proposed, where the leader is the authority of the energy system setting the prices with carbon emissions constraints and the followers are the energy users participating in the EIP.
The authors propose an iterative primal dual procedure based on KKT conditions, which guarantee global optimality upon convergence. The case studies show that this energy price system can positively affect both the economic benefit of the energy users and the environmental impact of the EIP. 

\section{Disaster prevention and response}
\label{sec:disaster}
Hazardous material (called “hazmat” hereafter) encompasses substances such as flam\-ma\-ble and explosive materials, posing significant human and environmental risks. Disaster prevention include hazmat transportation (see Table \ref{tab:hazmat}), tolls (see Table \ref{tab:haztoll}), and other governmental control policies (see Table  \ref{tab:hazcontrol}) introduced to address and mitigate these risks. Disaster response includes location, supply, and material distribution problems for emergency planning and relief logistics (see Table \ref{tab:emergency}).
These initiatives support SDG 1 (No Poverty), SDG 2 (Zero Hunger), and SDG 3 (Good Health and Well-being) by minimizing health hazards, mitigating the impact of disasters on communities and ensuring the efficient delivery of essential services including healthcare, food supplies, and emergency relief to affected areas.

\afterpage{
\clearpage
\newgeometry{margin=3cm}
\begin{landscape}
    \begin{table}[ht]
    \tiny
    \centering
    \caption{Summary of papers on bilevel applications focused on hazmat transportation. 
    }
    \label{tab:hazmat}
    \begin{tabularx}{\linewidth}{p{0.10\textwidth}p{0.19\textwidth}p{0.05\textwidth}p{0.02\textwidth}p{0.15\textwidth}p{0.15\textwidth}p{0.22\textwidth}p{0.18\textwidth}X}
    \toprule
    Application & References & \multicolumn{2}{c}{Bilevel game} & \multicolumn{2}{c}{Players} & \multicolumn{2}{c}{Objectives} & Method\\
    \cmidrule(lr){3-4}\cmidrule(lr){5-6}\cmidrule(lr){7-8}
    & & Class & Eq. & UL & LL & UL & LL &  \\
    \cmidrule(lr){1-1}\cmidrule(lr){2-2}\cmidrule(lr){3-4}\cmidrule(lr){5-6}\cmidrule(lr){7-8} \cmidrule(lr){9-9}       
    \multirow{10}{*}{\parbox{0.10\textwidth}{Hazmat transportation}} & \cite{kara2004} & SLMF  &  & Central authority & Carriers & Min total risk & Min distance & KKTR \\
    & \cite{GZARA201340} & SLMF  &  & Central authority & Carriers & Min total risk & Min distance & KKTR+CP \\
    & \cite{ERKUT20071389} & SLMF  &  & Central authority & Carriers & Min total risk & Min distance & GIA \\
    & \cite{ERKUT20082234} & SLMF  &  & Central authority & Carriers & Min total risk, Min total cost & Min distance & HEUR \\
    & \cite{BIANCO2009175} & SLSF  &  & Meta-local authority & Regional authority & Max risk equity & Min total risk & HEUR \\
    & \cite{xin2013} & SLMF  &  & Central authority & Carriers & Min total risk & Min distance & RO+HEUR \\
    & \cite{sun2016} & SLMF  &  & Central authority & Carriers & Min total risk & Min distance & RO+HEUR \\
    & \cite{liu2020} & SLMF  &  & Central authority & Carriers & Min setup cost, Min risk & Min transportation cost & CP+BD \\ 
    & \cite{TASLIMI2017} & SLMF  &  & Central authority & Carriers &  Min max risk & Min transportation cost & SLR,HEUR \\
    & \cite{Esfandeh2018} & SLMF &  & Central authority & Carriers &  Min risk & Min transportation cost &  HEUR \\
    \bottomrule
    \end{tabularx}
    \end{table}
    \vspace{-.6cm}
    \begin{table}[H]
    \tiny
    \centering
    \caption{Summary of papers on bilevel applications focused on hazmat toll setting. 
    }
    \label{tab:haztoll}
    \begin{tabularx}{\linewidth}{p{0.10\textwidth}p{0.19\textwidth}p{0.05\textwidth}p{0.02\textwidth}p{0.15\textwidth}p{0.15\textwidth}p{0.22\textwidth}p{0.18\textwidth}X}
    \toprule
    Application & References & \multicolumn{2}{c}{Bilevel game} & \multicolumn{2}{c}{Players} & \multicolumn{2}{c}{Objectives} & Method\\
    \cmidrule(lr){3-4}\cmidrule(lr){5-6}\cmidrule(lr){7-8}
    & & Class & Eq. & UL & LL & UL & LL &  \\
    \cmidrule(lr){1-1}\cmidrule(lr){2-2}\cmidrule(lr){3-4}\cmidrule(lr){5-6}\cmidrule(lr){7-8} \cmidrule(lr){9-9}       
    \multirow{3}{*}{\parbox{0.10\textwidth}{Hazmat toll setting}} & \cite{marcotte2009} & SLMF & & Government & Carriers & Min travel risk and cost & Max utility & Dual. SLR \\
    & \cite{ASSADIPOUR2016208} & SLMF & & Government & Carriers & Min travel risk, Min travel cost & Min cost & Hybrid Alg \\
    & \cite{LOPEZRAMOS2019} & SLMF & \checkmark & Network operator & Vehicles & Max profit & Min cost & Dual. SLR \\
    \bottomrule
    \end{tabularx}
    \end{table}
    \vspace{-.6cm}
    \begin{table}[H]
    \tiny
    \centering
    \caption{Summary of papers on bilevel applications focused on hazmat control policies. 
    }
    \label{tab:hazcontrol}
    \begin{tabularx}{\linewidth}{p{0.10\textwidth}p{0.19\textwidth}p{0.05\textwidth}p{0.02\textwidth}p{0.15\textwidth}p{0.15\textwidth}p{0.22\textwidth}p{0.18\textwidth}X}
    \toprule
    Application & References & \multicolumn{2}{c}{Bilevel game} & \multicolumn{2}{c}{Players} & \multicolumn{2}{c}{Objectives} & Method\\
    \cmidrule(lr){3-4}\cmidrule(lr){5-6}\cmidrule(lr){7-8}
    & & Class & Eq. & UL & LL & UL & LL &  \\
    \cmidrule(lr){1-1}\cmidrule(lr){2-2}\cmidrule(lr){3-4}\cmidrule(lr){5-6}\cmidrule(lr){7-8} \cmidrule(lr){9-9}
    \multirow{2}{*}{\parbox{0.10\textwidth}{Hazmat control policies}} 
    & \cite{CHIOU201616} & SLMF & \checkmark & Network operator & Vehicles & Min total travel delay & Min risk & CP \\
    & \cite{CHIOU2017} & SLMF & \checkmark & Network operator & Vehicles & Min total travel delay & Min risk  & CP \\
    & \cite{AMOUZEGAR1999100} & SLMF  &  & Central authority & Firms & Min total costs & Min cost & HEUR \\
    & \cite{BHAVSAR2022} & SLSF & & Government & Railroad operator & Min risk & Min cost & KKTR \\
    \bottomrule
    \end{tabularx}
    \end{table}
    \vspace{-.6cm}
    \begin{table}[H]
    \tiny
    \centering
    \caption{Summary of papers on bilevel applications focused on emergency planning and disaster response. 
    }
    \label{tab:emergency}
    \begin{tabularx}{\linewidth}{p{0.10\textwidth}p{0.19\textwidth}p{0.05\textwidth}p{0.02\textwidth}p{0.15\textwidth}p{0.15\textwidth}p{0.22\textwidth}p{0.18\textwidth}X}
    \toprule
    Application & References & \multicolumn{2}{c}{Bilevel game} & \multicolumn{2}{c}{Players} & \multicolumn{2}{c}{Objectives} & Method\\
    \cmidrule(lr){3-4}\cmidrule(lr){5-6}\cmidrule(lr){7-8}
    & & Class & Eq. & UL & LL & UL & LL &  \\
    \cmidrule(lr){1-1}\cmidrule(lr){2-2}\cmidrule(lr){3-4}\cmidrule(lr){5-6}\cmidrule(lr){7-8} \cmidrule(lr){9-9}
    \multirow{6}{*}{\parbox{0.10\textwidth}{Emergency planning and disaster response}} 
    & \cite{SAFAEI2018} & SLSF & & Relief SC designer & Relief SC operator & Min total relief SC cost & Min total supply risk &  Dual. SLR \\
    & \cite{li2019post} & SLSF &  & Emergency command center & Distribution centers administrator & Max road network accessibility & Max centers satisfaction; Min total delivery time & Hybrid GA \\ 
    & \cite{gao2022} & SLSF &  & Network designer  & Network operator & Min total dissatisfaction level & Min total expected time & SLR \\ 
    & \cite{liuluo2012}  & SLMF & \checkmark & Emergency manager  & Evacuees & Min total evacuation cost & Min perceived travel time & GA \\ 
    & \cite{YI2017}  & SLMF & \checkmark & Emergency manager & Residents & Min weighted travel time\&risk & Min travel time & HEUR \\ 
    & \cite{GUTJAHR2016} & SLMF & \checkmark & Aid-providing organization & Beneficiaries & Min opening cost; Min uncovered demand & Min weighted travel and unmet demand cost & $\epsilon$-const. + B\&B + FWA \\ 
    \bottomrule
    \end{tabularx}
\end{table}
\end{landscape}
\restoregeometry
}
\subsection{Hazmat transportation}
\label{subsec:haz-network}
In the domain of hazardous material transportation, road segment selection and network design are pivotal aspects. The main problem in the field is the Hazardous Materials Transportation Network Design Problem (HTNDP), which has been addressed with several bilevel programming approaches.
In the seminal work by \cite{kara2004}, the authors formulate the HTNDP under a bilevel framework where the leader is the government authority and the followers are the carriers.  
The leader determines road segments to be closed for carriers
while the followers make route choices.
The KKTR is linearized and solved by a commercial solver on a given instance. Results on a medium-sized case study in Western Ontario (Canada) show that a win-win situation exists for government and carriers, since both the risk and total distance decrease when a shipment-specific regulation for road segment selection is applied.
In \cite{GZARA201340}, the authors provide an exact CP method and a family of valid cuts for the bilevel HTNDP.
In \cite{ERKUT20071389}, the authors restrict the HTNDP by \cite{kara2004} to the case where only one route between any given origin and destination is allowed, resulting in a minimum risk hazmat tree network problem. 
They propose a greedy insertion algorithm (GIA) that iteratively adds paths to an optimal tree network, which is found by solving the SLR of the problem. The method only certifies the feasibility to the general problem. They test the GIA on a case study in Ravenna (Italy), consisting of a large network of around 100 nodes. 
They find different solutions with variable risk and cost, one of which reduces risk by 60\% with respect to the unregulated solution with only a 6\% increase in costs.
In \cite{ERKUT20082234}, the authors extend the HTNDP by \cite{kara2004} to the case of an undirected network and propose a bi-objective UL model. 
Results obtained with a heuristic using different risk measures on the case study of Ravenna as well as on random instances reveal that the heuristic is effective and time-efficient, both for the single-objective and the bi-objective problem. 
A different bilevel HTNDP, where the leader and the follower are, respectively, a local and a regional authority, is explored by \cite{BIANCO2009175}.
The study provides an iterative heuristic method based on \cite{ERKUT20082234}, which provides stable heuristic solutions when no optimistic/pessimistic approach is chosen. The authors demonstrate its effectiveness on a case study in Rome (Italy).

Addressing risk uncertainty is essential in hazmat transportation, where several external factors such as climate and road conditions can impact the hazmat risk. In the uncertain HTNDP, an interval of risk values is given for each arc of the network.
In \cite{xin2013}, intervals of risk values are introduced in the model proposed by \cite{kara2004} such that each commodity on each arc has its own set of possible risk values.  
The leader is the government and the followers are the carriers.
The authors provide a robust heuristic approach based on the minimax regret criterion and use it to find a robust shortest path for each commodity. The approach is tested on a small-size real-world case study in China for the transportation of solid and gas hazmat, showing good-quality solutions.
\cite{sun2016} present two robust bilevel models for the HTNDP, wherein the risk uncertainty associated with each arc is defined using intervals. A distinction is made between homogeneous and heterogeneous risk scenarios for all shipments.
The models incorporate an uncertainty budget, denoting the total number of arcs exposed to uncertainty for all shipments.
A heuristic algorithm based on a Lagrangian relaxation is proposed and tested on the Ravenna case study from  \cite{kara2004} and on larger instances from the city of Barcelona (Spain).
\cite{liu2020} combine the robust HTNDP and hazmat facility location under a bilevel framework, introducing uncertainty for hazmat risk and transportation demand. The leader decides which facilities to open and which road segments to close,
while the followers decide on their routes. 
The worst-case scenario is considered using an RO approach. The author adopt the CP method by \cite{GZARA201340}, combine it with a BD, and include uncertainty. The resulting exact method is proven to be effective compared to the direct solution of the SLR.
These studies collectively contribute to the understanding and optimization of the HTNDP, addressing environmental and social perspectives.
Other HTNDP extensions include, for example, equity risk (see \citealt{TASLIMI2017}) and time-dependent road closure policies (see \citealt{Esfandeh2018}).

\subsection{Hazmat toll setting}
\label{subsec:haz-toll}
Toll setting in hazmat transportation is an alternative tool for policy makers and road network operators to contain the hazmat risk exposure for road users.
\cite{marcotte2009} compare the bilevel HTNDP model by \cite{kara2004}, in which road segments are closed to hazmat carriers, with their toll setting model in which tolls are applied to roads for hazmat carriers. The government sets tolls on roads, while carriers decide their shipment routes. The study, which presents an SLR of the toll-setting model based on primal-dual optimality conditions, reveals that the toll setting model gives a lower optimal risk than the previous network design model.
\cite{ASSADIPOUR2016208} introduce a bi-objective bilevel model where the government deters the carriers from using certain rail intermodal terminals by assigning a toll to each hazmat container passing through them, and the carriers select the routes. 
The authors address the problem with a hybrid algorithm where a GA for multi-objective optimization is combined with the optimal solution of the LL problem. 
A case study in USA is solved by comparing the toll setting model with the network design model (where the government closes some terminals instead of setting tolls). Different results suggest that a two-stage procedure (toll setting first and network design after) obtains promising solutions.
Network design and toll setting are sometimes integrated into a unique bilevel setting, such as in the model by \cite{LOPEZRAMOS2019}. The authors formulate a strong-duality-based SLR, which is then heuristically solved on benchmark instances.

\subsection{Hazmat control policies}
\label{subsec:haz-pol}
Other bilevel applications for hazmat transportation deal with control policies. 
\cite{CHIOU201616} introduces nonlinear bilevel models for hazmat transportation with signal-controlled road networks. The aim is to regulate hazmat traffic by setting signals. The LL problem involves a traffic equilibrium. 
The author proposes an iterative CP method based on bundling subgradients from previous iterations.  \cite{CHIOU2017} extends the problem to incorporate travel demand uncertainty.
Addressing pollution control policies, \cite{AMOUZEGAR1999100} present a bilevel model for hazardous waste management (see Section \ref{par:haz-waste}).
Subsidies are another effective tool for hazmat transportation risk control. \cite{BHAVSAR2022} introduce an SLSF model, where the government offers subsidies to induce the railroad operator to take alternative routes that are away from high-risk network links. 
The authors present a KKTR, applied to a real case in USA, demonstrating that even modest subsidies can result in significant risk reduction.

\subsection{Emergency planning and disaster response}
\label{subsec-emergency}
Natural disasters require efficient relief supply chains to provide relief commodities to communities, including food, clothing, and medicines, among many.
Leader and follower of SLSF bilevel games model the strategic and operational levels of decisions of relief logistics decision makers in various contexts.
\cite{SAFAEI2018} propose a robust model to optimize the supply/demand process of relief commodities under uncertainty. In the UL, the leader decides the locations of transshipment transfer depots close to the disaster areas for collecting commodities from central warehouse and transfer decisions to the disaster areas.
In the LL, the transfer of commodities from suppliers to the central warehouse are determined.
The authors obtain a KKTR. Results on flood disaster real scenarios in Iran demonstrate the effectiveness of the proposed method. 
\cite{li2019post} formulate the multi-period bilevel road network repair work scheduling and relief logistics problem for earthquake disaster relief. The leader makes the restoration strategy (i.e., repair crew assignment and routing) for each period, whereas the follower 
makes relief logistics and relief material delivery decisions according to the UL strategy.
The proposed hybrid GA is able to solve a case study in China in a short computational time,  with a resulting 15\% level of road repair at the end of the time horizon. 
The multi-commodity rebalancing problem allows to rebalance surplus and shortage of relief commodities over the transportation network to satisfy the potential demand at all relief centers. 
\cite{gao2022} presents a stochastic bilevel model under demand and transportation network availability uncertainty. In the UL, the incoming and outgoing shipments at relief centers are determined to achieve fairness between the centers, while in the LL, the routing decisions are made.
The authors solve an SLR on an earthquake case study in China, showing applicable and effective decision-making results.

Prevention or post-disaster decisions of the policy maker and the consequent behaviour of individuals affected by a disaster are well modeled by SLMF games with equilibrium between the followers.
\cite{liuluo2012} study the emergency evacuation process, in which signalized and uninterrupted flow intersections in the evacuation network must be optimally located for traffic crossing-elimination and signal control. 
In the UL, the emergency manager decides the intersections location. In the LL problem, modeled as a SUE problem, evacuees make evacuation routes decisions. The authors propose a GA with the successive weighted averages method for the convex LL equilibrium and solve a case study in China, obtaining total evacuation time reductions of up to 40\%. 
\cite{YI2017} present a bilevel model with a similar structure (SLMF with dynamic traffic equilibrium in the LL), where the leader decides when and where to issue orders and the followers in the LL are the residents deciding if, when, and how to evacuate.
The UL is a multi-stage stochastic program with hurricane occurrence uncertainty. The authors propose a heuristic algorithm based on a Lagrangian relaxation with scenario decomposition and solve a case study in North Carolina for different hurricane scenarios.
\cite{GUTJAHR2016} study the FLP of distribution centers for a natural disaster response humanitarian SC. 
The leader makes location opening decisions and the followers make supply decisions in a traffic user equilibrium model.
The authors propose an exact approach based on the $\epsilon$-constraint method for the bi-objective UL, a B\&B, and the Frank-Wolfe algorithm (FWA) for the LL Wardrop equilibrium. Medium-size instances for a case study of rural communities in Senegal are solved to optimality, showing reductions of up to 40\% of the unmet demand.  

\section{Conclusions and future research directions}
\label{sec:concl}
Bilevel optimization models hierarchical decision-making processes involving multiple actors.
The field has been rapidly expanding in both solution methodologies and applications, approaching different real-world situations involving policy makers, industrial entities, and end users. 
In alignment with contemporary public and private concerns surrounding the ``3P'' framework, that stands for ``profit, planet, people'', and the well-known Sustainable Development Goals (SDGs), this paper analyzes the literature on bilevel applications with sustainability perspectives, considering the three-dimensional aspect of the term.
By delving into the literature, we identify bilevel applications that incorporate at least two (out of three) dimensions of sustainability – economic, environmental, and social – within the upper or lower levels of the hierarchical game structure.
These applications align with SDGs by promoting resilient infrastructure, sustainable industrialization, and low-carbon technologies, addressing urban sustainability challenges including disaster risk reduction, air pollution control, and efficient transportation networks, and encouraging waste reduction and environmentally conscious production practices in businesses. 

We study the players of these bilevel games, their decisions and goals, along with the proposed solution methods, and the resulting insights on real-world case studies for transportation and logistics, production planning and manufacturing, waste, water, and agriculture management, supply chains, and disaster prevention and response.
Applications cover renowned mathematical programming problems, such as assignment, location, and routing problems, among many, which have only recently been extended to incorporate optimization of environmental and social aspects.

In bilevel applications, leaders typically emerge as authorities designing service systems for the private/public domain or companies optimizing their businesses to provide consumers with products. Optimization problems encompass diverse decision-making processes, ranging from carbon policies to green vehicle and route selection, waste reduction, environmental preservation, clean energy production, health risk mitigation from hazardous materials, and emergency logistics optimization.
A recurring pattern observed in our analysis is that the upper-level problem, typically led by an authority, incorporates green and social goals alongside economic ones. Conversely, the lower-level problem of end users usually prioritizes efficiency or profitability. Numerous case studies covering industrial and public environments across various countries in Europe, America, Asia, and Africa showcase the potential of bilevel optimization in achieving sustainability goals. Positive outcomes include reduced carbon emissions, expanded and efficient network infrastructures, and maximized social benefits measured in terms of people's preferences and health risk mitigation.

Most of these bilevel applications adopt the optimistic approach or remain unspecified when using heuristic methods. 
Pessimistic models, where the leader decides against the worst-case scenario, are still rare in sustainability applications.
Most existing pessimistic models focus on interdiction or risk-averse investment problems, which, in the contexts covered by our review, typically do not relate to sustainability actions.

For single-leader-single-follower applications, a single-level reformulation can often be derived, allowing the problem to be solved using commercial solvers (when possible) or through tailored heuristics on larger instances. When multiple objectives are present  (typically at the upper level), evolutionary algorithms are commonly employed to obtain heuristic bilevel solutions. 
Similarly, single-level reformulations are sometimes achievable for single-leader-multi-follower problems; however, evolutionary and genetic algorithms are widely used, particularly when multiple objectives exist at one or both levels. 
For problems involving a lower-level equilibrium among the followers, genetic algorithms and tailored heuristics are currently the most effective methods, although some attempts at single-level reformulations have also been explored.

Despite the wealth of applications, there remain aspects that can be further investigated. We sketch possible directions for future research in the field of bilevel optimization with sustainability perspectives:
\begin{enumerate}
    \item Future research should focus on developing frameworks that integrate environmental and social perspectives into the lower level of bilevel optimization models,  addressing issues such as environmental awareness, community well-being, social equity, and the broader societal impacts of business decisions. This would allow to investigate the role of end users, including companies, consumers, and citizens in shaping sustainable practices, acknowledging not only external influences like public policies but also the individual motivations and values that drive decision-making.
    For instance, a bilevel approach could be applied in the fleet management context, where authorities and companies, respectively, with emission regulations and green practices reduce the impact of vehicular emissions but also improve the economic efficiency of the fleet.
    \item The exploration of the social dimension in bilevel optimization remains in its early stages, as evidenced by the fact that SDGs 4 (Quality Education), 5 (Gender Equality), 10 (Reduced Inequalities), and 16 (Peace, Justice, and Strong Institutions) are not addressed by any of the papers included in our review. Researchers should investigate how decisions, especially of the leader of the bilevel game, impact various social aspects, including equity, accessibility, and well-being considerations, as well as stakeholder engagement and corporate responsibility. This involves the incorporation of more social indicators and metrics into bilevel goals. 
    The work by \cite{HuWQLG22} is an exception in the field of bilevel facility location, tackling the location and size problem for general service infrastructures from a social perspective. Their bilevel model incorporates proximity to desirable and undesirable areas to evaluate as a metric of location environmental impacts, and optimizes service efficiency and density as measures of user satisfaction.
    \item Fairness in transportation has been neglected so far in bilevel optimization applications for transportation and logistics. Interesting future research should investigate the modeling of fairness-related requirements and objective functions to consider leader's social optimum that is also fair in terms of, for instance, the difference between the maximum and minimum travel time of followers. 
    \item  Given the variety of routing problems studied in the literature, the field of bilevel green routing can be further explored. The few works reported in this survey open room for future works on hierarchical routing decisions, such as the introduction of green policies and vehicle selection for routing plans with reduced carbon emissions.
    \item Other sectors, such as healthcare and medicine, should be explored with bilevel optimization. Very few applications have emerged, for instance, for home healthcare supply chain optimization and genetic studies on the human body network structures.
    \item From the methodological perspective, we notice that the majority of sustainable problems with a bilevel structure and formulation are still solved using heuristics, and certain lax practices persist in the field (see \citealt{CAMACHOVALLEJO2024}). 
    In many heuristic-based approaches, the focus is on obtaining high-quality solutions for the follower’s problem rather than solving it to optimality. In that case, there is no guarantee that the obtained solution is even bilevel-feasible, which is a significant drawback of metaheuristic approaches. 
    In addition, the treatment of multiple follower optimal solutions is neglected, leading to potentially ill-posed bilevel formulations. Therefore, the gap between the rigorous analysis of bilevel mathematical programs and the use of efficient metaheuristics for practical problem-solving is currently, to the best of our knowledge, an open question in bilevel literature.
    \item We expect that recent advances in computational bilevel optimization (see, e.g., \citealt{KleinertLLS21}) and future algorithmic developments will allow to tackle these complex bilevel settings using exact algorithms as well. In particular, comparison of approaches based on single-level reformulations, which leverage optimality conditions, provides an interesting direction for future research and a reference in the literature. Scalability is another crucial factor (in addition to optimality) to be considered in the future development of exact solution techniques.  Addressing the challenges posed by large-scale bilevel problems remains a key priority, as scalability is a crucial factor in the development of exact solution techniques.
\end{enumerate}
In summary, bilevel optimization has already made substantial contributions to sustainable decision-making, and we expect that this collection of works not only consolidates existing knowledge but also lays the foundation for future developments in the research field and supports practitioners in more structured and holistic decisions.

Another important sector that, due to the large number of possible applications, has remained out of the scope of this study, is that of electricity and energy markets. The literature on this sector is already large and well structured, including cases that adopt the pessimistic approach (see, e.g., \citealt{ALVES2018}). Our future research will be dedicated to the analysis of green bilevel energy applications.

{\footnotesize
\noindent \textbf{Disclosure of interest} The authors report there are no competing interests to declare. \\
\noindent \textbf{Acknowledgements} Manuel Iori and Giulia Caselli gratefully acknowledge financial support under grant PRIN2022PNRR – M4C2INV1.1, NextGenerationEU - Award 1409/2022 -  Project Calipso. The authors also thank the three anonymous referees for their valuable comments, which helped improve this work.
}

\appendix
\section{Acronyms}
\label{app:acro}
\begin{table}[H]
    \centering
    \caption{Acronyms}
    \label{tab:acronyms}
    \resizebox{0.75\textwidth}{!}{ 
    \begin{tabular}{llll}
        \toprule
         Acronym & Meaning & Acronym & Meaning  \\
         \midrule
AON & All-or-nothing	&	LP &  Linear programming 	\\
B\&B & Branch-and-Bound	&	MILP & Mixed-integer linear programming	\\
B\&C & Branch-and-Cut	&	MINLP & Mixed-integer nonlinear programming	\\
BD & Benders decomposition	&	MLMF & Multi-leader-multi-follower	\\
BSA & Binary search algorithm	&	MLSF & Multi-leader-single-follower	\\
CCG & Column-and-constraint generation	&	MO & Multi-objective	\\
CG & Column generation	&	MOEA & Multi-objective evolutionary algorithm	\\
CP & Cutting plane	&	MOPSO &  Multi-objective particle swarm optimization	\\
DP &  Dynamic programming	&	NSGA-II & Non-dominated Sorting Genetic Algorithm II	\\
EA & Evolutionary algorithm	&	OD & Origin-destination	\\
EIP & Eco-industrial park	&	PRP & Pollution routing problem	\\
EV & Electric vehicle	&	PSO &  Particle swarm optimization	\\
FLP & Facility location problem	&	RO &  Robust optimization	\\
FWA & Frank-Wolfe algorithm	&	SA & Simulated annealing	\\
GA & Genetic algorithm	&	SC & Supply chain	\\
GHG & Greenhouse gas	&	SCND & Supply chain network design	\\
GIA & Greedy insertion algorithm	&	SDG & Sustainable development goal	\\
GNEP & Generalized nash equilibrium problem	&	SLMF & Signle-leader-multi-follower	\\
GWP & Global warming potential	&	SLR  &  Single-level reformulation	\\
HEUR &  (Meta)heuristic	&	SLSF & Single-leader-single-follower	\\
HTNDP & Hazardous materials transportation network design problem	&	SO & Stochastic optimization	\\
IA & Iterative algorithm	&	SUE & Stochastic user equilibrium	\\
ILS &  Iterated local search	&	TOD & Transit-oriented development	\\
KKT &  Karush–Kuhn–Tucker 	&	TS & Tabu search 	\\
KKTR &  Karush–Kuhn–Tucker reformulation	&	UL & Upper level	\\
LL & Lower level 	&	VRP & Vehicle routing problem	\\
LLVFR & Lower-level value function reformulation	&	WS & Weighted sum	\\
\bottomrule
    \end{tabular}
}
\end{table}

\singlespacing
\bibliographystyle{elsarticle-harv} 
\small
\bibliography{cas-refs}






\end{document}